\documentclass[reqno]{amsart}
\usepackage{color}
\usepackage{amssymb,amsmath,amsthm,amstext,amsfonts}
\usepackage[dvips]{graphicx}
\usepackage{psfrag}
\usepackage{url}

\makeatletter \@addtoreset{equation}{section} \makeatother

\renewcommand\thetable{\thesection.\@arabic\c@table}

\newtheorem{theorem}{Theorem}[section]
\newtheorem{proposition}{Proposition}[section]
\newtheorem{lemma}{Lemma}[section]
\newtheorem{corollary}{Corollary}[section]
\newtheorem{definition}{Definition}[section]
\newtheorem{remark}{Remark}[section]
\newtheorem{example}{Example}[section]

\newcommand{\htop}{h_{\topp}}

\newcommand{\al} {\alpha}       
\newcommand{\be} {\beta}        
\newcommand{\ga} {\gamma}    \newcommand{\Ga}{\Gamma}
\newcommand{\de} {\delta}       \newcommand{\De}{\Delta}
\newcommand{\vep}{\varepsilon}

\newcommand{\la} {\lambda}      \newcommand{\La}{\Lambda}

\newcommand{\si} {\sigma}

\newcommand{\R}{\mathbb{R}}

\newcommand{\diam}{\operatorname{diam}}
\newcommand{\dist}{\operatorname{dist}}

\newcommand{\topp}{\operatorname{top}}

\newcommand{\Leb}{Leb}
\newcommand{\ti}{\tilde }
\newcommand{\Ptop}{P_{\topp}}

\newcommand{\mL}{{\mathcal L}}

\newcommand{\cC}{\mathcal{C}}

\newcommand{\cP}{\mathcal{P}}
\newcommand{\cL}{\mathcal{L}}

\newcommand{\cM}{\mathcal{M}}

\newcommand{\cF}{\mathcal{F}}

\def\proc#1{\medbreak\noindent{\it #1}\hspace{1ex}\ignorespaces}
\def\ep{\noindent{\hfill $\Box$}}

\begin{document}

\title{Non-uniform specification and large deviations for weak Gibbs measures}

\author{PAULO VARANDAS}

\address{Departamento de Matem\'atica, Universidade Federal da Bahia\\
  Av. Ademar de Barros s/n, 40170-110 Salvador, Brazil.}
\email{paulo.varandas@ufba.br\\ 
\url{www.pgmat.ufba.br/varandas} }

\date{\today}

\begin{abstract}
We establish bounds for the measure of deviation sets associated to
continuous observables with respect to not necessarily invariant weak 
Gibbs measures. Under some mild assumptions, we obtain upper and lower bounds for the
measure of deviation sets of some non-uniformly expanding maps,
including quadratic maps and robust multidimensional non-uniformly
expanding local diffeomorphisms.
For that purpose, a measure theoretical weak form of specification is 
introduced and proved to hold for the robust classes of multidimensional
nonuniformly expanding local diffeomorphisms and Viana maps. 
\end{abstract}


\maketitle

%%%%%%%%%%%%%%%%%%%%%%%%%%%%%%%%%%%%%%%%%%%%%%%%%%
\section{Introduction }

The theory of Large Deviations concerns the study of the rates of
convergence at which time averages of a given sequence of random
variables converge to the limit distribution. An application of
these ideas into the realm of Dynamical Systems is useful to
estimate the velocity at which typical points of ergodic invariant
measures converge to the corresponding space averages. More
generally, given a continuous transformation $f$ on a compact metric
space $M$ and a reference measure $\nu$, one would like to provide
sharp estimates for the $\nu$-measure of the deviation sets
$$
\left
\{
    x \in M : \frac{1}{n} \sum_{j=0}^{n-1} g(f^j(x)) > c
\right\}
$$
for all continuous functions $g: M \to \R$ and real numbers $c$. To
this purpose, a priori estimates on the measure of the dynamical
balls
$$
B(x,n,\vep)
    =\left\{y\in M : d(f^j(y),f^j(x)) \leq \vep, \;\forall \, 0\leq j \leq n\right\}
$$
for $x \in M$, $\vep>0$ and $n\geq 1$ are useful and somewhat
necessary since points that belong to the same dynamical ball have
nearby Birkhoff averages with respect to continuous functions.

Some large deviations ideas and techniques are particularly useful
to the study of the thermodynamical formalism of transformations
with some hyperbolicity. Recall that the variational principle for
the pressure asserts that for every continuous potential $\phi$
$$
\Ptop(f,\phi)=\sup \left\{ h_\eta(f)+\int \phi \;d\eta \right\},
$$
where the supremum is taken over all invariant probability measures
$\eta$. A measure $\mu$ that attains the supremum in the variational
principle is called an \b{equilibrium state} for $f$ with respect
to the potential $\phi$. A large deviations theory was developed for
uniformly hyperbolic systems restricted to a basic piece of the
non-wandering set and H\"older continuous potentials in both
discrete and time-continuous settings. Indeed,
such hyperbolic transformations 
admit a unique equilibrium state with respect to
any H\"older continuous potential (see~\cite{Si72,Bo75,Ru76b}),
Young, Kifer and Newhouse~\cite{You90,Ki90,KN91}) established, in
the mid nineties, large deviation principles for this important open
class of dynamical systems: the rate of decay is given explicitly in
terms of the distance of all invariant measures $\eta$ with ``bad"
space averages to the equilibrium state $\mu$. Focusing on the
discrete time case, the sharp lower and upper bounds obtained in
~\cite{You90} for the measure of deviation sets yield as a
consequence that for any ergodic equilibrium state $\mu$ and every
continuous observable $g$, the measure of the set of points whose
time average $\frac{1}{n} \sum_{j=0}^{n-1} g(f^j(x))$ is far from
the space average $\int g \,d\mu$ decreases exponentially fast.  Two
key ingredients to obtain the large deviations principle are that
equilibrium states are Gibbs measures and that,
when restricted to a basic piece of the non-wandering set, every
uniformly hyperbolic dynamical system is semi-conjugated to a
subshift of finite type that satisfies a very ``strong mixing"
condition known as specification. This notion, introduced by
Bowen~\cite{Bo71}, means roughly that any finite sequence of pieces
of orbit can be well approximated by periodic ones.

Our purpose here is to give a contribution for the ergodic theory of non-uniformly expanding 
maps in two directions. Namely, we introduce a measure-theoretical non-uniform specification property 
and obtain upper and lower large deviations bounds with respect to weak Gibbs measures as we now detail.

On the one hand,
the specification property constitutes an important tool in dynamical systems which is useful e.g. 
to obtain uniqueness of equilibrium states for expansive transformations, to study large deviations or to study the 
multifractal formalism for associated to Birkhoff averages.   However, and despite the fact that the later property holds for topologically mixing interval maps and dynamical systems with arbitrary small finite Markov partition,  conceptually one cannot expect this to hold with  great generality in the absence of uniform hyperbolicity. 
For that reason we introduce a measure theoretical non-uniform specification property and prove that
it holds for a large class of robust nonuniformly expanding maps as in \cite{VV10} and the multidimensional nonuniformly hyperbolic attractors with critical region considered in \cite{Vi97}. Such class of transformations that 
may not satisfy the strong specification property seem to constitute first multidimensional examples  presenting a
weak form of specification in a nonuniformly hyperbolic context. One should mention that other mild forms of specification were introduced by Saussol, Troubetzkoy, Vaienti~\cite{STV03} to study the relation between recurrence and dimension in dynamical systems, by Pfister, Sullivan \cite{PS05} and Thompson~\cite{Th10}
to the study of multifractal formalism for Birkhoff averages associated to beta-shifts, and also by 
Yamamoto~\cite{Yam09} to study large deviations for automorphisms on compact metric abelian groups.

On the other hand, since the nineties many efforts have been made in the attempt to
extend the theory of large deviations to the scope of non-uniformly
hyperbolic dynamics and some important results in that direction
have been obtained recently. Ara\'ujo and Pac\'ifico \cite{AP06}
established large deviation upper bounds for the deviation sets of
physical measures for non-uniformly expanding maps (in the sense of
\cite{ABV00}). More recently, Melbourne and Nicol \cite{MN08,Mel09}
studied systems that admit some inducing Markov structure, and
proved that the measure of points with atypical time averages for a
H\"older continuous potential has the same decay rate as the
inducing time itself. In particular, less than exponential rate of
convergence to equilibrium is studied. Independently, in the case of
exponential tail, Rey-Bellet and Young~\cite{RY08} obtained similar
and sharper results. The construction of (countable) expanding
Markov maps in \cite{Pi09} provides many examples where the previous
results apply.  Large deviations principles were also obtained by Yuri  (see~\cite{Yu05})
in the context of shifts with countably many symbols and by 
Comman and Rivera (see ~\cite{CoRi08}) for non-uniformly expanding rational maps.
Notice that in \cite{AP06} the authors establish large deviation upper bounds
with respect to Lebesgue measure, while in \cite{MN08,Mel09,RY08,Yu05,CoRi08} the decay rate
for the measure of deviation sets is studied with respect to the invariant probability measure.
More recently, Chung~\cite{Chu11} obtained some large deviation principles for the Lebesgue measure
on Markov maps satisfying some technical conditions of similar flavor to our nonuniform specification property.

Hence, to the best of our knowledge, the theory of large deviations with respect to not necessarily 
invariant reference measures arising from thermodynamical formalism is far from complete.  

Inspired by the pioneering work of Young~\cite{You90} our purpose in this direction is to obtain large deviations 
estimates for non-uniformly expanding maps that exhibit the non-uniform specification property with respect to a not necessarily invariant weak Gibbs measure.
Weak Gibbs measures are such that  the measure of dynamical balls $B(x,n,\vep)$ is given according to 
the $n$th Birkhoff sums of some potential on the orbit of $x$ up to some multiplicative constant which has at 
most subexponential growth in $n$. See \eqref{eq.weakgibbs} below for the precise definition.
Moreover, such measures arise naturally in the thermodynamical formalism of many non-uniformly expanding 
maps, where equilibrium states arise as invariant measures absolutely continuous with respect to some
reference weak Gibbs measure as illustrated in Section~\ref{s.applications}.
Roughly, one proves that the set of points whose time averages remain far from the space average with respect to
the equilibrium measure decrease exponentially fast, with a decay rate which is related
to the existence of invariant expanding measures with frequent hyperbolicity (we refer the reader to
Theorem~\ref{thm.weakGibbs.weak} for the precise statement).  
Since equilibrium states associated to uniformly expanding dynamics and H\"older continuous
potentials satisfy the strong Gibbs property our results partially extend the ones of~\cite{You90} 
to the non-uniformly expanding setting. In particular, the thermodynamical approach used
here fails in the same extent to provide sharp large deviation bounds if there is non-uniqueness of equilibrium
states. See for instance \cite{Ki90} for an example in the uniformly hyperbolic setting. 

One should also mention that the large deviation results presented here are indeed complementary and extend 
the ones of obtained by Ara\'ujo and Pac\'ifico \cite{AP06} in the non-uniformly expanding setting. First, the reference measure is not necessarily Lebesgue and our assumptions in our Theorem~\ref{thm.weakGibbs.weak} rely on the 
Gibbs property for the reference measure while the assumptions of \cite{AP06} rely on the non-uniform hyperbolicity
and slow recurrence condition to the critical region. Finally, we obtain large deviation lower bound estimates which
were not available even in the case of the Lebesgue measure.

This paper is organized as follows. Our main results are stated
along Section~\ref{s.statements}. In Section~\ref{s.preliminaries}
we recall some necessary definitions and prove some preliminary
lemmas. The proofs of our main results are given in
Sections~\ref{s.proof1} and ~\ref{s.proof2}. 
Finally, in Section~\ref{s.applications} we present some examples 
and further questions.

%%%%%%%%%%%%%%%%%%%%%%%%%%%%%%%%%%%%%%%%%%%%%%%
\section{Statement of results }\label{s.statements}

%%%%%%%%%%%%%%%%%%%%%%%%%%%%%%%%%%%%%%%%%%%%%%%%
\subsection{Abstract Theorem}

Let $f:M\to M$ be a continuous transformation on a compact metric
space $M$ and let $\nu$ be some (not necessarily invariant)
probability measure. In this section we state an abstract result on
the deviation of Birkhoff averages given by continuous observables.

Given an observable $\phi: M \to \mathbb R$, we denote by
$S_n\phi(x)=\sum_{j=0}^{n-1}\phi\circ f^j$ the $n$th Birkhoff sum of
$\phi$. Given a full $\nu$-measure set $\La\subset M$, denote by
$\cF(\La)$ the set of continuous functions $\psi \in C(M,\R)$ so
that,
there exists $\de_0>0$ and for every $x \in \La$ and $0<\vep<\de_0$
there exists a sequence of positive constants $(K_n)_{n\geq1}$ 
such that $\lim\limits_{n\to\infty} \frac{1}{n} \log K_n(x,\vep)=0$
and
\begin{equation}\label{eq:abstract.measure.control}
K_n(x,\vep)^{-1} \;e^{-S_n\psi(x)}
    \leq \nu(B(y,n,\vep)) \leq
K_n(x,\vep) \;e^{-S_n\psi(x)}
\end{equation}
for every $n \geq 1$ and every $y \in M$ satisfying $B(y,n,\vep)
\subset B(x,n,\de_0)$. This is a generalization of the usual notion of Gibbs measure corresponding
which can be obtained e.g. in the case that $\Lambda$ is compact and $x\mapsto K_n(x,\vep)$
is continuous and independent of $n$. Here we do not assume the compactness of $\Lambda$
nor any regularity of the functions $K_n$. We define also
$\de(\vep,\beta)$ as the exponential decay rate corresponding to the $\nu$-measure of
the points whose constants $K_n$ grow at most $\beta$-exponentially,
that is, if
\begin{equation}\label{eq.set.Delta.n}
\De_{n}(\beta) 
    =\Big\{ x\in \La : K_n(x,\vep)<e^{\beta n}\Big\},
\end{equation}
then set $\de(\vep,\beta)=\limsup_{n\to\infty} \frac1n \log \nu(\De_n^c(\beta))$.
For notational simplicity, when no confusion is possible we shall omit the dependence on $\beta$ 
in the definition of the set $\De_n$. In a context of nonuniform hyperbolicity
the quantity $\de(\vep,\beta)$ appears as  the exponential decay of the
instants of hyperbolicity, does not depend on $\vep$ and it is negative
for interesting class of examples that appear in Section~\ref{s.applications}. 
Finally, the \emph{relative entropy} of an $f$-invariant probability measure $\eta$ 
is defined as $h_{\nu}(f,\eta)=\eta\mbox{-esssup}
\;h_{\nu}(f,\cdot)$, where
\[
h_{\nu}(f,x)
    =\lim_{\vep \to 0} \limsup_{{\substack{n \to \infty}}}-\frac{1}{n}\log \nu(B(x,n,\vep)),
    \quad\text{for all}\; x \in M.
\]
We will also need the following:

\begin{definition}
\label{d.strong.specification}
We say that a map $f$ satisfies the \emph{specification property} if
for any $\vep>0$ there exists an integer $N=N(\vep)\geq 1$ such that
the following holds: for every $k\geq 1$, any points $x_1,\dots,
x_k$, and any sequence of positive integers $n_1, \dots, n_k$ and
$p_1, \dots, p_k$ with $p_i \geq N(\vep)$ 
there exists a point $x$ in $M$ such that
$$
\begin{array}{cc}
d\Big(f^j(x),f^j(x_1)\Big) \leq \vep, &\forall \,0\leq j \leq n_1
\end{array}
$$
and
$$
\begin{array}{cc}
d\Big(f^{j+n_1+p_1+\dots +n_{i-1}+p_{i-1}}(x) \;,\; f^j(x_i)\Big)
        \leq \vep &
\end{array}
$$
for every $2\leq i\leq k$ and $0\leq j\leq n_i$.

\end{definition}

Note that this notion of specification is purely topological and
is slightly weaker than the one introduced by Bowen~\cite{Bo71}, 
that requires that any finite sequence of pieces of orbit is well approximated by periodic orbits.
In fact, this condition is known to imply that the system is topologically mixing \cite{Bo71}. 
It might seem that specification is quite rare among most dynamical systems. However,
Blokh \cite{Bl83} proved in a surprising way that the notions of
specification and topologically mixing coincide for every
one-dimensional \emph{continuous} mapping. This is no longer true if
the one-dimensional map fails to be continuous (see
e.g.~\cite{Buz97}). We refer the reader to \cite{Wa82} for more details on the specification
property. Our first result is as follows.

\begin{theorem}\label{thm.weakGibbs.strong}
Assume that $\htop(f)<\infty$, let $\nu$ be a probability measure
and let $\La\subset M$ be such that $\nu(\La)=1$. Given $g \in
C(M,\R)$ and $c \in \R$, if $\psi \in \cF(\La)$ then for every small
$\vep,\beta>0$ it holds
$$
 \limsup_{n \to \infty} \frac{1}{n}\log \nu\bigg[x \in M : \frac{1}{n} S_n g(x) \geq c\bigg]
    \leq \max\left\{\de(\vep,\beta)\; , \sup \big\{h_\eta(f) - \int \psi \,d\eta\big\}\!+\!\beta\right\}
$$
where the supremum is over all invariant
probability measures $\eta$ such that $\int g \,d\eta \geq c$.
Moreover, it holds that
$$
\liminf_{n \to \infty} \frac{1}{n}\log \nu \left(x \in M :
\frac{1}{n} S_n g(x) > c\right) \geq \sup \Big\{ h_\eta(f) -
h_\nu(f,\eta) \Big\} ,
$$
where the supremum is taken over all
\emph{ergodic} measures $\eta$ satisfying $\int g \,d\eta > c$.
Furthermore if $\psi \in \cF(\La)$ and $f$ satisfies the
specification property then 
$$
 \liminf_{n \to \infty}
    \frac{1}{n}\log \nu\left( x\in M :\frac{1}{n} S_n g(x) > c\right)
    \geq \sup \left\{h_\eta(f) - \int \psi\,d\eta\right\},
$$
where the supremum is taken over all invariant
probability measures $\eta$ such that $\eta(\La)=1$
and $\int g \,d\eta
> c$.
\end{theorem}

This theorem 
generalizes Theorem~1 in \cite{You90}, where $\La=M$ was assumed 
to be compact and some uniform control on the measure of partition elements was required.

%%%%%%%%%%%%%%%%%%%%%%%%%%%%%%%%%%%%%%%%%%%%%%%%
\subsection{Deviation bounds for non-uniformly expanding maps}

%%%%%%%%%%%%%%%%%%%%%%%%%%%%%%%%%%%%%%%%%%%%%%%%
\subsubsection{Context} 

Let $M$ be a compact Riemaniann manifold and let $f : M \to M$ be a
$C^{1+\al}$ local diffeomorphism outside of a compact 
critical or singular region $\cC$. Assume:
\begin{itemize}
\item[(H)] $f$ {\em behaves like a power of the distance to
            the critical or singular set $\cC$}:  there exist
            $B>1$ and $\be\in (0,1)$ such that for every $x,y\in M\setminus \cC$ with
            $\dist(x,y)<\dist(x,\cC)/2$ and every $v\in T_x M$:
            \vspace{.1cm}
            \begin{itemize}
            \item[(a)]  $\frac{1}{B}\dist(x,\cC)^{\be}\leq \frac{\|Df(x)v\|}{\|v\|}\leq
                        B\dist(x,\cC)^{-\be};$
                        \vspace{.1cm}
            \item[(b)] $\left|\log\|Df(x)^{-1}\|- \log\|Df(y)^{-1}\|\:\right|\leq
                        B\frac{\dist(x,y)}{\dist(x,\cC)^{\be}};$
                        \vspace{.1cm}
            \item[(c)] $\left|\log|\det Df(x)^{-1}|-\log|\det Df(y)^{-1}|\:\right|\leq B
                        \frac{\dist(x,y)}{\dist(x,\cC)^{\be}}.$
                        \vspace{.1cm}
            \end{itemize}
\end{itemize}
This condition was proposed in \cite{ABV00} as a multidimensional
counterpart of the non-flat critical points in one-dimensional
dynamics.  We also assume the following condition on $f$:

\vspace{.1cm} \noindent
(C) There exists $L>0$ and $\ga \in (0,1)$ 
such that for any small $\vep>0$ every connected component in the preimage of a set of diameter $\vep$ is 
contained in a ball of radius $L \vep^{\ga}$.
\vspace{.1cm}

This condition is clearly satisfied if $f$ is a local diffeomorphism and, since 
$f$ behaves like a power of the distance to $\cC$,  it is most likely to hold e.g. if
$\cC$ has empty interior.  
Such condition is satisfied by the class transformations with singularities (quadratic  
and Viana maps) considered in Section~\ref{s.applications}.
 Let $\phi:M\backslash \cC \to \R$ be a H\"older continuous potential and assume:

\begin{enumerate}
\item[(P1)] There exists a probability measure $\nu$ that is positive on open
        sets, it is non-singular with respect to $f$ with H\"older continuous
        Jacobian $J_\nu f=\la e^{-\phi}$, for some $\la>0$.
        We will refer to $\nu$ as a \emph{conformal measure} associated to $\phi$;
\item[(P2)] $(f,\nu)$ has \emph{non-uniform expansion}: there exists $\si>1$ such
        that for $\nu$-a.e. $x$ 
        \begin{equation*}
        \limsup_{n\to \infty} \frac{1}{n} \sum_{j=0}^{n-1} \log\|
        Df(f^j(x)^{-1}\| \leq -2\log \si<0
        \end{equation*}
        and
        \begin{equation*}
        (\forall \vep>0)\, (\exists \de>0) \; 
        \limsup_{n\to \infty} \frac{1}{n} \sum_{j=0}^{n-1} -\log \dist_\de(f^j(x),\cC) <\vep,
        \end{equation*}
\end{enumerate}
where for any given $\delta>0$, we let $\dist_\delta(x,\cC)$ be
the $\delta$-truncated distance from a point $x$ to $\cC$ defined as
$\dist (x,\cC)$ if $\dist(x,\cC)< \delta$ and equal to $1$ otherwise.

These assumptions are quite natural in a context of non-uniform
hyperbolicity and are verified by a large class of maps and
potentials. For instance, if $f$ is a non-uniformly expanding map
(in the sense of \cite{ABV00}) and $\phi=-\log |\det Df|$ then the
Lebesgue measure is a conformal measure that satisfies (P1) and
(P2). Usually conformal measures appear as eigenmeasures associated
to the dual $\cL^*_\phi$ of the Ruelle-Perron-Frobenius operator
\begin{equation*}\label{eq. RPF operator}
\mL_\phi g(x) = \sum_{f(y)=x} e^{\phi(y)} \, g(y),
\end{equation*}
acting on the space of probability measures $\cM$. Moreover,
hypothesis (P1) and (P2) together with the fact that the potential
$\phi$ is H\"older continuous yield that $\nu$ is a \emph{weak Gibbs
measure}: there are $P \in \R$ and $\de>0$ so that for any
$0<\vep\leq\de$ and almost every $x$ there is a sequence of positive
numbers $(K_n)_{n \geq 1}$ (depending also on $\phi$) satisfying
$\lim\limits_{n\to\infty} \frac{1}{n} \log K_n(x,\vep)=0$ and
\begin{equation}\label{eq.weakgibbs}
K_n(x,\vep)^{-1}
    \leq \frac{\nu(B(x,n,\vep))}{e^{-P n + S_{n}\phi(y)}}
    \leq K_n(x,\vep)
\end{equation}
for every $y \in B(x,n,\vep)$ (see e.g. \cite{VV10}). Compare to
Lemma~\ref{le:measure.control} and
Corollary~\ref{co:measure.control} below. We say that $n$ is a
\textit{$(\si,\de)$-hyperbolic time} for $x\in M$ (or hyperbolic
time for short) if there is a small positive constant $b>0$ such
that
\begin{equation*}
\prod_{j=n-k}^{n-1} \| Df(f^j(x)^{-1}\| \leq \si^{-k} \quad
\text{and} \quad \dist_\de(f^{n-k}(x),\cC)>\si^{-bk}
\end{equation*}
for every $1 \leq k \leq n$. The non-uniform expansion condition
(P2) guarantees the existence of infinitely many hyperbolic times
$\nu$-almost everywhere. We refer the reader to Subsection~\ref{Sec.
Hyperbolic Times} for more details. Let $H$ denote the set of points
with infinitely many hyperbolic times, $n_1(\cdot)$ be the first
hyperbolic time map and $\Gamma_n=\{x\in M: n_1(x)>n\}$.
We say that a probability measure $\eta$ is \emph{expanding} if
$\eta(H)=1$. In particular, any invariant expanding measure has only
positive Lyapunov exponents.
We also assume:

\vspace{.1cm}
\begin{enumerate}
\item[(P3)] There is a unique equilibrium state $\mu$ for $f$
        with respect to $\phi$, it is absolutely continuous with respect to $\nu$,
         there exists a positive constant $K>0$ such that the density satisfies 
        $d\mu/d\nu \ge K^{-1}$, and $n_1 \in L^1(\mu)$.
\end{enumerate}
\vspace{.1cm}

The last assumption above essentially means that the decay of the
first hyperbolic time map is at least polynomial of order
$n^{-(1+\vep)}$, for some $\vep>0$. We refer the reader to the works
\cite{ABV00,BS03,Yu03,Yu05,OV08,VV10}, just to quote some classes of
maps and potentials that satisfy our assumptions.

%%%%%%%%%%%%%%%%%%%%%%%%%%%%%%%%%%%%%%%%%%%%%%%%
\subsubsection{Non-uniform specification property}

In contrast to the topological concept of specification we introduce a
{\bf measure theoretical} notion. 

\begin{definition}
We say that $(f,\mu)$ satisfy the {\bf non-uniform specification property}
if there exists $\de>0$ such that for $\mu$-almost every $x$ and
every $0<\vep<\de$ there exists an integer $p(x,n,\vep)\geq 1$
satisfying
$$
\lim_{\vep\to 0}\limsup_{n \to \infty} \frac{1}{n} p(x,n,\vep)=0
$$
and so that the following holds: given points $x_1, \dots, x_k$ in a full $\mu$-measure set
and positive integers $n_1, \dots, n_k$, if $p_i \geq p(x_i,n_i,\vep)$
then there exists $z$ that $\vep$-shadows the orbits of each $x_i$
during $n_i$ iterates with a time lag of $p(x_i,n_i,\vep)$ in
between $f^{n_i}(x_i)$ and $x_{i+1}$, that is,
$$
z \in B(x_1,n_1,\vep)
    \quad\text{and}\quad
f^{n_1+p_1+\dots +n_{i-1}+p_{i-1}}(z) \in B(x_i,n_i,\vep)
$$
for every $2\leq i\leq k.$

\end{definition}

These notions means 
that almost every finite pieces of orbits are approximated by a real
orbit such that the time lag between two consecutive pieces of
orbits is small proportion of the size of the piece of orbit being
shadowed. 
Clearly, if the strong specification property holds then $(f,\eta)$
satisfies the non-uniform specification property for every $f$-invariant
probability measure $\eta$. 
Let us also mention that a notion of non-uniform specification 
property similar to the one introduced in \cite{STV03} (using Pesin theory)
would also be enough to obtain the lower bound estimates in Theorem
~\ref{thm.weakGibbs.weak} below. We shall not use or prove this fact here. 

In opposition to the specification property we expect this weak form
of specification to hold in a broad non-uniformly hyperbolic setting. 
We refer the reader to Section~\ref{s.applications} for some examples
in which the later condition holds but may fail to satisfy the specification property.

%%%%%%%%%%%%%%%%%%%%%%%%%%%%%%%%%%%%%%%%%%%%%
\subsubsection{Deviation bounds for non-uniformly expanding maps}

The following result extends the large deviation results proven in  \cite{You90} for 
uniformly hyperbolic maps.

\begin{theorem}\label{thm.weakGibbs.weak}
Let $M$ be a compact manifold and $f:M \to M$ be a $C^{1+\al}$ local
diffeomorphism outside a critical or singular region $\cC$ that
satisfies (H) and (C). Let $\phi: M\setminus \cC \to \R$ be an H\"older
continuous potential and let $\nu$ and $\mu$ be probability measures
given by (P1)-(P3). If $g \in C(M,\R)$ and
$c \in \R$ then 
\begin{align*}
 \limsup_{n \to \infty} \frac{1}{n} & \log \nu\left(x \in M :\frac{1}{n} S_n g(x) \geq c\right)\\
        & \leq \max
        \left\{
            \sup \left\{-P+h_\eta(f) + \int \phi \,d\eta\right\},
            \limsup_{n\to\infty} \frac1n \log \mu(\Gamma_n)
        \right\}
\end{align*}
where the supremum is taken over all invariant
probability measures $\eta$ such that $\int g \,d\eta \geq c$. If,
in addition, $f$ satisfies the specification property or $(f,\mu)$ satisfies the
non-uniform specification property then
$$
 \liminf_{n \to \infty} \frac{1}{n}\log \nu\left( x\in M :\frac{1}{n} S_n g(x) > c\right)
    \geq \sup \left\{-P+h_\eta(f) + \int \phi\,d\eta\right\},
$$
where the supremum is taken over all invariant
probability measures $\eta$ such that
$\eta(H)=1$, $\int g \,d\eta > c$ and
$n_1 \in L^1(\eta)$.
\end{theorem}

In consequence one can estimate the decay of the deviation set as follows:

\begin{corollary}\label{cor:conseq}
Under the previous assumptions,
\begin{align*}
 \limsup_{n\to \infty} \frac1n & \log \nu \left(x \in M : \left|\frac1n S_n g(x)-\!\int g\,d\mu\right| \ge c\right)\\
        & \leq \max
        \left\{
            \sup \left\{-P+h_\eta(f) + \int \phi \,d\eta\right\},
            \limsup_{n\to\infty} \frac1n \log \mu(\Gamma_n)
        \right\}
\end{align*}
where the supremum is taken over all invariant
probability measures $\eta$ such that $|\int g \,d\eta -\int g d\mu|\geq c$,
and
\begin{align*}
\liminf_{n\to \infty} \frac1n & \log \nu \left(x \in M : \left|\frac1n S_n g(x)-\!\int g\,d\mu\right|>c\right)\\
        &\geq \sup \left\{-P+h_\eta(f) + \int \phi\,d\eta\right\}
\end{align*}
where the supremum is taken over all invariant
probability measures $\eta$ such that
$\eta(H)=1$, $|\int g \,d\eta -\int g d\mu|>  c$ and
$n_1 \in L^1(\eta)$.
\end{corollary}

Some comments are in order.
First notice that the upper bound estimate takes into account the loss of uniform expansion in terms of the
decay of the first hyperbolic time map. In particular, if the first hyperbolic time map fails to have
exponential decay then the right hand side in the previous upper bound is zero,
since the other term is also non-positive. In \cite{MN08} less than
exponential deviations are proven for systems that admit a Young tower with 
inducing time has polynomial decay. More recently, in \cite{AFLV10} a relation between the rate of decay of correlations and the large deviations with respect to the invariant measure has been established. This
reenforces the idea that a condition on the tail of the first hyperbolic time map should not be easily removed in general. See Example~\ref{ex.maneville} for a more detailed discussion.
Moreover, we expect Theorem~\ref{thm.weakGibbs.strong} 
to hold in the more general setting of zooming measures introduced in \cite{Pi09} 
since our ingredients are bounded distortion and growth to large scale. 
Finally, these results should also extend to the nonuniformly hyperbolic setting, e.g. 
the class of partially hyperbolic diffeomorphisms with contracting direction and center-unstable 
mostly expanding direction introduced in~\cite{ABV00}. Since the construction of general
equilibrium states for such class of maps is still not available, the lack of motivating examples
lead us to state the results only in the non-uniformly expanding setting.

%%%%%%%%%%%%%%%%%%%%%%%%%%%%%%%%%%%%%%%%%%%%%%
\section{Preliminary results}\label{s.preliminaries}

%%%%%%%%%%%%%%%%%%%%%%%%%%%%%%%%%%%%%%%%%%%%%%%
\subsection{Metric Entropy}\label{s.entropy}

First we recall some definitions. Let $\vep>0$ and $n\geq 1$ be
arbitrary. A set $E\subset M$ is \emph{$(n,\vep)$-separated} if
$d_n(x,y)>\vep$ for every $x, y \in E$ with $x\neq y$, where $d_n:M
\times M \to \R_0^+$ is the metric given by
$$
d_n(x,y)
 =\max\limits_{0\leq j\leq n-1} d(f^j(x),f^j(y)).
$$
If, in addition, $E$ has maximal cardinality we say that it is a
\emph{maximal $(n,\vep)$-separated} set. Note that for any maximal
$(n,\vep)$-separated set $E$, the dynamical balls $B(x,n,\vep)$
centered at points in $E$ are pairwise disjoint and that the union
$\cup_{x\in E} B(x,n,2\vep)$ covers $M$. We recall some properties
of topological and metric entropy.

\begin{proposition}\cite{Bo71}
Let $f:M\to M$ be a continuous map in a compact metric space $M$. If
$\htop(f)$ denotes the topological entropy of $f$ then
$$
\htop(f)= \lim_{\vep \to 0} \limsup_{n\to \infty}
        \frac1n \log N(n,\vep),
$$
where $N(n,\vep)$ the minimum number of $(n,\vep)$ dynamical balls
necessary to cover $M$.
\end{proposition}

A metric counterpart of this result is as follows. Let $\eta$ be an
invariant probability measure and $\de>0$ arbitrary. Given $\vep>0$
let $N(n,\vep,\de)$ be the minimum number of $(n,\vep)$-dynamical
balls necessary to cover a set of measure larger than $1-\de$.

\begin{proposition}\cite[Theorem I.I]{Ka80}\label{p.Katok}
Let $f:M\to M$ be a homeomorphism in a compact metric space $M$ and
$\eta$ an $f$-invariant probability measure. Hence
$$
h_\eta(f)= \lim_{\vep \to 0} \limsup_{n\to \infty}
        \frac1n \log N(n,\vep,\de)
        = \lim_{\vep \to 0} \liminf_{n\to \infty}
        \frac1n \log N(n,\vep,\de),
$$
for every $\de>0$.
\end{proposition}
\vspace{.2cm}

%%%%%%%%%%%%%%%%%%%%%%%%%%%%%%%%%%%%%%%%%%%
\subsection{Hyperbolic times }\label{Sec. Hyperbolic Times}

In this subsection we recall some properties of hyperbolic times.

\begin{definition}
We say that $(f,\eta)$ is \emph{non-uniformly expanding} if there
exists $N \geq 1$ and $\si>1$ such that almost every $x$ satisfies
\begin{equation}\label{nue1}
\limsup_{n\to \infty} \frac{1}{n} \sum_{j=0}^{n-1} \log\|
Df^N(f^{jN}(x)^{-1}\| \leq -2\log \si<0
\end{equation}
and the \emph{slow recurrence condition}: for every $\vep>0$ there
exists $\de>0$ such that for $\mu$-almost every point $x \in M$ it
holds that
\begin{equation}\label{nue2}
 \limsup_{n \to \infty} \frac{1}{n} \sum_{j=0}^{n-1} -\log \dist_\de(f^j(x), \cC) < \vep.
\end{equation}
\end{definition}

Let $B,\beta$ be given by condition $(H2)$ and take
$0<b<\{\frac{1}{2}, \frac{1}{2\beta}\}$.
A sufficiency criterium for the existence of hyperbolic times is
given as application of Pliss' lemma.

\begin{lemma}\cite[Lemma 5.4]{ABV00} \label{positive density}
There exists constants $\theta>0$ and $\de>0$ (depending only on $f$
and $c$) such that if $x \in M\backslash \cup_n f^n(\cC)$ satisfies
(\ref{nue1}) and (\ref{nue2}) then the following holds: for every
large $N\geq 1$ there exist a sequence of integers $1 \leq n_1(x) <
n_2(x) < \dots< n_l(x) \leq N$, with $l \geq \theta n$ so that
\begin{equation}\label{eq. ht}
\prod_{j=n-k}^{n-1} \| Df(f^j(x)^{-1}\| \leq \si^{-k}
    \quad
\text{and} \quad \dist_\de(f^{n-k}(x),\cC)>\si^{bk}.
\end{equation}
\end{lemma}

One of the main features of hyperbolic times is stated below.

\begin{lemma}\cite[Lemma 2.7]{ABV00}\label{delta}
Given $c>0$ and $\de>0$ there exists a constant
$\delta_1=\delta_1(c,\de,f)>0$ such that if $n$ is a hyperbolic time
for a point $x$ then $f^n$ maps diffeomorphically the dynamical ball
$V_n(x)=B(x,n,\delta_1)$ onto the ball $B(f^n(x),\delta_1)$ around $f^n(x)$
and radius $\delta_1$ and
\begin{equation*}\label{eq. backward distances contraction}
d (f^{n-j}(y),f^{n-j}(z))
        \leq \si^{-\frac{j}{2}}\;  d(f^n(y),f^n(z))
\end{equation*}
for every $1\leq j \leq n$ and every $y,z \in V_n(x)$.
\end{lemma}

Using that $J_\nu f=\la e^{-\phi}$ is H\"older continuous and the
backward distances contraction at hyperbolic times we obtain a
bounded distortion property.

\begin{corollary}\label{c.bounded.distortion}
There exists $K_0>0$ such that for every $y,z \in V_n(x)$
$$
K_0^{-1}
    \leq \frac{J_\nu f^n(y)}{J_\nu f^n(z)}\leq
K_0.
$$
\end{corollary}

%%%%%%%%%%%%%%%%%%%%%%%%%%%%%%%%%%%%%%%%%%%%%%%
\subsection{Control of the measure of dynamical balls}

Now we prove a useful lemma on the measure of dynamical balls for
weak Gibbs measures. 
In what follows $\de_1$ stands for the diameter of the hyperbolic
ball as in Lemma~\ref{delta}.

\begin{lemma}\label{le:measure.control}
For every $0<\vep<\de_1$ there exists a positive constant
$K(\vep)>0$ such that if $n$ is a hyperbolic time for $x$ and
$B(y,n,\vep) \subset B(x,n,\de)$ then
    $$
    K(\vep)^{-1}
    \leq \frac{\nu(B(y,n,\vep))}{e^{-Pn+S_n\phi(y)}} \leq
    K(\vep),
    $$
where $P=\log\la$.
\end{lemma}

\proc{Proof}
One has $f^n(B(y,n,\vep))=B(f^n(y),\vep)$ by backward distance
contraction at hyperbolic times. Hence,
Corollary~\ref{c.bounded.distortion} asserts that
$$
1 \ge \nu(B(f^n(y),\vep))
        = \int_{B(y,n,\vep)} e^{-S_n\psi} \, d\nu
        \ge K_0^{-1} e^{Pn-S_n\phi(y)} \, \nu (B(y,n,\vep)).
$$
Using that $\nu$ is positive on open sets and the compactness of $M$
it follows that the measure of every ball of radius $\vep$ is
bounded away from zero. Thus the other inequality is obtained analogously.
\ep\medbreak

The following very interesting consequence is that dynamical balls
have comparable measure that only depends on the center of the ball.

\begin{corollary}\label{co:measure.control}
Assume that $x\in H$. For every $0<\vep<\de_1$ and $n \ge 1$ there
exists a positive constant $K_n(x,\vep)>0$ such that if $B(y,n,\vep)
\subset B(x,n,\de)$ then
    $$
    K_n(x,\vep)^{-1}
    \leq \frac{\nu(B(y,n,\vep))}{e^{-Pn+S_n\phi(y)}} \leq
    K_n(x,\vep).
    $$
\end{corollary}

\proc{Proof}
Given an arbitrary $n$ write $n_i(x) \leq n < n_{i+1}(x)$, where
$n_i$ and $n_{i+1}$ are consecutive hyperbolic times for $x$. Using
that $B(y,n,\vep) \subset B(y,n_i(x),\vep)$ it is clear that
\begin{align*}
\nu(B(y,n,\vep))
        & \le K(\vep) e^{(\sup|\phi|+|P|)(n-n_i(x))} \, e^{-Pn+S_{n}\phi(y)} \\
        & \le K_n(x,\vep) \, e^{-Pn+S_{n}\phi(y)},
\end{align*}
with $K_n(x,\vep)=K(\vep)\exp[(\sup|\phi|+|P|) (n-n_i(x))]$ (depends
only on the center $x$). This finishes the proof of the corollary.
\ep\medbreak

Now we prove that the constants $K_n$ have subexponential growth with respect to
every invariant expanding measure such that the
first hyperbolic time map is integrable.  More precisely,

\begin{lemma}\label{le:subexpgrowth}
Let $\eta$ be an $f$-invariant and expanding probability measure so
that $n_1 \in L^1(\eta)$ and let $K_n(x,\vep)$ be given as above. Then,
\begin{equation}\label{eq:subexpgrowth}
\lim\limits_{\vep\to 0} \limsup\limits_{n \to \infty}
    \frac1n \log K_n(x,\vep)=0
\end{equation}
for $\eta$-almost every $x$. In consequence, $\psi=\phi-P$ belongs
to $\cF(H)$.
\end{lemma}

\proc{Proof}
This proof resembles the one of Proposition~3.8 in~\cite{OV08}. Let
$\eta$ be an $f$-invariant, expanding probability measure so that
$n_1 \in L^1(\eta)$ and take $\beta>0$ arbitrary. Given $x\in H$, $n
\ge 1$ and $0<\vep<\de_1$ recall that $K_n(x,\vep)\le K(\vep)
\exp[(|P|+\sup|\phi|) n_1(f^{n_i(x)}(x))]$, where $n_i(x) \le n \le
n_{i+1}(x)$ are consecutive hyperbolic times for $x$. Set
$C_\beta=\beta n /(|P|+\sup|\phi|)-\log K(\vep)$. If
$K_n(x,\vep)>e^{\beta n}$ this implies that $n_1(f^{k}(x))> C_\beta
n \geq C_\beta \,k$, where $k=n_i(x)$. This shows that
\begin{align*}
\{ x\in H : K_n(x,\vep)>e^{\beta n} \; \text{i.o.}\}
    & \subset \{ x\in H : n_1(f^n(x),\vep)>e^{\beta n} \;\text{i.o.}\}\\
    & \subset \bigcup_{n \ge 1} \{x \in H : n_1(f^{n}(x))> C_\beta n \}
\end{align*}
Furthermore, using the invariance of $\eta$ and the integrability
assumption
\begin{align*}
\sum_{n=1}^{+\infty} \eta\big(x \in H : n_1(f^n(x))> C_\beta n\big)
        =\sum_{n=1}^{+\infty} \eta\big(x\in H : n_1(x)> C_\beta n\big)
         \leq \int n_1 \, d\eta
         <\infty.
\end{align*}
Using Borel-Cantelli lemma this proves that $K_n(x,\vep)\le e^{\beta
n}$ for all but finitely many values of $n$ for $\eta$-almost every
$x$. Since $\beta$ was taken arbitrary, this completes the proof of
the first claim above.

Using that $n_1 \in L^1(\mu)$ and $d\mu/d\nu$ is bounded from below
by a constant (recall assumption (P3)) it follows that
\eqref{eq:subexpgrowth} holds $\nu$-almost everywhere. Together with
Corollary~\ref{co:measure.control} this shows that $\psi=\phi-P$
belongs to $\cF(H)$ and that $\nu$ is a weak Gibbs measure. This
finishes the proof of the lemma.
\ep\medbreak

\begin{remark}
\label{rmk:hyperbolic.times}
It follows from \eqref{eq:subexpgrowth} and the definition of the
constants $K_n(x,\vep)$ that if $n_1 \in L^1(\eta)$ then given
$\beta>0$, for $\eta$-almost every $x$ there exists $n_x\ge 1$ such
that $n-n_i(x) \le \beta n$ for every $n \ge n_x$. In fact we prove
even more: given $\beta>0$ then for $\eta$-almost every $x$ there
exists $n_x\ge 1$ such that $n_{i+1}(x)-n_i(x) \le \beta n$ for
every $n \ge n_x$.
\end{remark}

%%%%%%%%%%%%%%%%%%%%%%%%%%%%%%%%%%%%%%%%%%%%%
\section{Abstract deviation bounds}\label{s.proof1}

In this section we prove Theorem~\ref{thm.weakGibbs.strong}. Upper and lower 
bounds for the measure of the deviation sets are given separately.

%%%%%%%%%%%%%%%%%%%%%%%%%%%%%%%%%%%%%%%%%%%%%%
\subsection{Upper bound}

Let $g \in C(M,\R)$, $c\in\R$ and $\psi \in \cF(\La)$ be fixed. We
want to prove that for every small $\vep,\beta>0$
$$
 \limsup_{n \to \infty} \frac{1}{n}\log \nu\left[x \in M : \frac{1}{n} S_n g(x) \geq c\right]
    \leq \max\left\{\de(\vep,\beta)\; , \sup \big\{h_\eta(f) - \int \psi \,d\eta\big\}\!+\!\beta\right\}
$$
where the supremum is taken over all invariant probability measures
$\eta$ such that $\int g \,d\eta \geq c$. We use the following
result from Calculus (see e.g. \cite[Lemma 9.9]{Wa82}).
\begin{lemma}
Given $n\geq 1$, real numbers $(a_i)_{i=1..n}$ and $0\leq p_i\leq 1$
such that $\sum_{i=1}^n p_i=1$ then
$$
\sum_{i=1}^n p_i(a_i-\log p_i)
        \leq \log \left(\sum_{i=1}^n e^{a_i} \right),
$$
and the equality holds if and only if $p_i=\frac{e^{a_i}}{\sum_j
e^{a_j}}$.
\end{lemma}

\noindent Let $B_n$ denote the set of points $x\in M$ so that $S_n
g(x) \geq c n$. Recall that $\La$ is a $\nu$-full measure set and,
for every $x\in \La$ and every small $\vep>0$ it holds that
$$
\nu(B(x,n,\vep))\leq K_n(x,\vep) e^{-S_n\psi(x)},
$$
with $\limsup_{n} \frac1n \log K_n(x,\vep)=0$. Let $\beta>0$ and
$0<\vep<\de_0$ be arbitrary small and $n \ge 1$ be fixed. Then $B_n
\subset \De_n^c \cup \left(B_n \cap \De_n\right)$, where $\De_n$ is
as in~\eqref{eq.set.Delta.n}. Moreover, if $E_n\subset B_n \cap
\De_{n}$ is a maximal $(n,\vep)$-separated set, $B_n\cap \De_{n}$ is
contained in the union of the dynamical balls $B(x,n,2\vep)$
centered at points of $E_n$ and, consequently,
\begin{equation}\label{eq.upper}
\nu(B_n) < 
        \nu(\De_n^c) + e^{\beta n} \sum_{x \in E_n} e^{-S_n\psi(x)}
\end{equation}
for every $n$. Now, consider the probability measures $\si_n$ and
$\eta_n$ given by
$$
\si_n
   =\frac{1}{Z_n} \sum_{x \in E_n} e^{-S_n\psi(x)} \de_x
   \quad\text{and}\quad
\eta_n
    = \frac{1}{n} \sum_{j=0}^{n-1} f^j_* \si_n,
$$
where $Z_n=\sum_{x \in E_n} e^{-S_n\psi(x)}$, and let $\eta$ be an
weak$^*$ accumulation point of the sequence $(\eta_n)_{n}$. It is
not hard to check that $\eta$ is an $f$-invariant probability
measure. Assume $\cP$ is a partition of $M$ with diameter smaller
than $\vep$ and $\eta(\partial\cP)=0$. Each element of $\cP^{(n)}$
contains at most one point of $E_n$. By the previous lemma
$$
H_{\si_n}(\cP^{(n)}) - \int S_n \psi \,d\si_n
    = \log \Big( \sum_{x\in E_n} e^{-S_n\psi(x)} \Big)
$$
which, as in the usual proof of the variational principle (see
\cite[Pages 219-221]{Wa82}), guarantees that
\begin{equation}\label{eq.upper2}
\limsup_{n \to \infty}
    \frac{1}{n} \log Z_n
    \leq h_\eta(f) - \int \psi \, d\eta.
\end{equation}
Observe also that $\int \psi\,d\eta\geq c$ by weak$^*$ convergence
since $E_n$ is contained in $B_n$ and
$$
 \int g \,d\eta_n
        = \frac{1}{n} \sum_{j=0}^{n-1} \frac{1}{Z_n} \sum_{x \in E_n} e^{S_n\phi(x)}.\,g\circ f^j(x)
        \geq c.
$$
Finally, it follows from \eqref{eq.upper} and \eqref{eq.upper2} that
for every $\beta>0$
\begin{align*}
\limsup_{n \to \infty} \frac{1}{n} \log \nu(B_n)
   & \leq \max\left\{\de(\vep,\beta)\; ,  h_\eta(f) - \int \psi \, d\eta+\beta\right\}\\
   & \leq \max\left\{\de(\vep,\beta)\; ,  \sup\left\{h_\xi(f) - \int \psi \, d\xi\right\}+\beta\right\},
\end{align*}
where the supremum is over all invariant probability measures. This
completes the proof of the first statement in
Theorem~\ref{thm.weakGibbs.strong}.

%%%%%%%%%%%%%%%%%%%%%%%%%%%%%%%%%%%%%%%%%%%%%%
\subsection{Lower bound using ergodic measures}\label{ss.lowerA}

Let $g:M\to \R$ be a continuous map, take $c\in \R$ and a 
$\beta>0$ small. If $\eta$ is an ergodic probability measure 
such that $\int g \,d\eta > c$ we claim that
\begin{equation*}
\liminf_{n \to \infty}
        \frac{1}{n}\log \nu \left(x \in M : \frac{1}{n} S_n g(x) > c\right)
        \geq
        h_\eta(f) - h_\nu(f,\eta)-2\beta.
\end{equation*}
Denote by $B_n$ the set of points $x\in M$ such that $S_n g(x) > cn$
and fix $\de_2=\frac12 (\int g d\eta-c)$. Notice that $h_\eta(f)\leq
\htop(f)<\infty$ and that we may assume $h_\nu(f,\eta)<\infty$
(because otherwise there is nothing to prove). Hence $\eta$-almost
every point $x$ satisfies
\begin{equation}
h_\nu(f,x)
    =\lim_{\vep\to 0} \limsup_{n\to \infty} -\frac1n \log \nu(B(x,n,\vep))
    \leq h_\nu(f,\eta)<\infty.
\end{equation}
Since $\eta$ is ergodic then $\frac{1}{n} S_n g(x) \to \int
g\,d\eta$ for $\eta$-almost every $x$. Choose $\xi>0$ by uniform
continuity so that $|g(x)-g(y)|<\de_2$ whenever $d(x,y)<\xi$.
Observe that if $n_0 =n_0(\beta)\geq 1$ is large and $\de\in(0,\xi)$
is small enough then the set $D$ of points $x\in M$ satisfying
\begin{equation}\label{eq.1}
\frac{1}{n} S_n g(x) > c+\de_2
    \quad\text{and}\quad
\nu (B(x,n,\vep)) \geq e^{-[h_\nu(f,\eta)+ \beta]n}
\end{equation}
has $\eta$-measure at least $\frac12$, and that the minimal number
$N(n,2\vep,\frac12)$ of $(n,2\vep)$-dynamical balls necessary to
cover a set of $\eta$-measure at least $\frac12$ satisfies
\begin{equation}\label{eq.3}
N\left(n,2\vep,\frac12\right)
    \geq e^{[h_\eta(f)- \beta]n}
\end{equation}
for every $n\geq n_0$ and every $0<\vep\leq \de$. Indeed, the
existence of such $n_0$ and $\de$ is a consequence of
Proposition~\ref{p.Katok}, the definition of relative entropy and
ergodicity. Moreover, it follows from our choice of $\xi$ and
$\de_2$ that
\[
B_n \supset
    \bigcup_{x \in  D} B(x,n,\vep)
    \supset  D
\]
for all $n\ge n_0$ and $0<\vep<\xi$. So, if $\vep>0$ is small and
$E_n \subset D$ is a maximal $(n,\vep)$-separated set, using that
the dynamical balls $B(x,n,\vep)$ centered at points in $E_n$ are
pairwise disjoint contained in $B_n$ and the union $\cup_{x\in E_n}
B(x,n,2\vep)$ covers $D$, relations~\eqref{eq.1} and \eqref{eq.3}
yield that 
$$
\nu(B_n)
    \geq \nu\Big( \bigcup_{x\in E_n} B(x,n,\vep) \Big)
    \geq \sum_{x\in E_n} \nu\Big( B(x,n,\vep) \Big)
    \geq e^{[h_\eta(f) - h_\nu(f,\eta) - 2\beta ]n}
$$
whenever $n \geq n_0$, which proves our claim. The second second assertion
in Theorem~\ref{thm.weakGibbs.strong} follows from the arbitrariness of $\beta$.

\begin{remark}
Since $B_n \supset B_n \cap \De_n$, where $\De_n$ is as
in~\eqref{eq.set.Delta.n} then
 $
 \nu(B_n) \ge \nu(B_n \cap \De_n).
 $
However, we have no estimate whatsoever for the measure of the
intersection in terms of $\nu(\De_n)$. Hence the previous result
shows only that the measure of the points with predetermined
Birkhoff averages decreases at most exponentially fast.
\end{remark}

%%%%%%%%%%%%%%%%%%%%%%%%%%%%%%%%%%%%%%%%%
\subsection{Lower bound over all invariant measures}

The proof of the last statement in
Theorem~\ref{thm.weakGibbs.strong} is divided in two steps. First we
prove the lower bound when the supremum is restricted over ergodic
measures. Afterwards we deduce the general bound using that every
invariant measure can be approximated by a finite collection of
ergodic measures and the specification to ``glue" together finite
pieces of orbits. We begin by the following lemma.

\begin{lemma}\label{l.lower.ergodic}
If $g\in C(M,\R)$, $c\in \R$ and $\psi \in \cF(\La)$ then
$$
 \liminf_{n \to \infty}
    \frac{1}{n}\log \nu\left( x\in M :\frac{1}{n} S_n g(x) > c\right)
    \geq h_\eta(f) - \int \psi\,d\eta
$$
for every ergodic probability measure $\eta$ such that $\eta(\La)=1$
and $\int g \,d\eta > c$.
\end{lemma}

\proc{Proof}
Fix $g\in C(M,\R)$, $c\in \R$ and $\psi \in \cF(\La)$, and denote by
$B_n$ the set of points $x\in M$ such that $S_n g(x)>cn$. Let
$\beta>0$ be a small constant and $\de_2=\frac12 (\int g d\eta-c)$.
Let $\xi>0$ be given by uniform continuity such that
$|g(x)-g(y)|<\de_2$ for any points $x,y\in M$ at distance smaller
than $\xi$. As before, if $n_0$ is large enough and $0<\vep<\xi$ is
small then the set $D$ of points $x\in \Lambda$ satisfying
\begin{equation}\label{eqq.1}
\frac{1}{n} S_n g(x) > c+\de_2,
    \quad
K_n(x,\vep)^{-1} \geq e^{-\beta n}
    \quad\text{and}\quad
\frac 1 n S_n \psi(x) < \int \psi \,d\eta+\beta
\end{equation}
for every every $n \geq n_0$ has $\eta$-measure at least $\frac12$
and $N\left(n,2\vep,\frac12\right) \geq e^{(h_\eta(f)- \beta)n}$ for
every $n\ge n_0$. Then, using that
 $
B_n \supset
    \bigcup_{x \in  D} B(x,n,\vep)
    \supset  D
 $
it follows that
$$
\nu(B_n)
    \geq \sum_{x\in E_n} \nu\Big( B(x,n,\vep) \Big)
    \geq e^{ [  h_\eta(f)-\int \psi \,d\eta- 3\beta ] n }
$$
for every maximal $(n,\vep)$-separated set $E_n \subset D$. This
proves that
$$
 \liminf_{n \to \infty}
    \frac{1}{n}\log \nu\left( x\in M :\frac{1}{n} S_n g(x) > c\right)
    \geq h_\eta(f) - \int \psi\,d\eta-3\beta.
$$
Since $\beta$ was taken arbitrary the statement in the lemma follows
directly.
\ep\medbreak

The following result asserts that any invariant probability measure
can be approximated by a finite convex combination of ergodic
measures supported in $\Lambda$. 

\begin{lemma}\label{l.ergodic.approximation}
Let $\eta=\int \eta_x d\eta(x)$ be the ergodic decomposition an
$f$-invariant probability measure $\eta$ such that $\eta(\La)=1$.
Given $\beta>0$ and a finite set $(\psi_j)_{1 \leq j\leq r} \subset
C(M)$ of continuous functions, there are positive real numbers
$(a_i)_{1 \leq i \leq k}$ satisfying $a_i\leq 1$ and $\sum a_i=1$,
and finitely many points $x_1, \dots, x_k$ such that the ergodic
measures $\eta_i=\eta_{x_i}$ from the ergodic decomposition satisfy
\begin{itemize}
\item[(i)] $\eta_i(\La)=1$;
\item[(ii)] $h_{\hat\eta}(f) \geq h_{\eta}(f) - \beta$; and
\item[(iii)] $|\int \psi_j \;d\hat\eta - \int \psi_j \;d\eta|
            <\beta$ for every $1 \leq j \leq r$;
\end{itemize}
where $\hat \eta= \sum\limits_{i=1}^k a_i \eta_i$.
\end{lemma}

\proc{Proof}
Fix the $f$-invariant probability measure $\eta$ such that
$\eta(\La)=1$. By ergodic decomposition theorem and convexity of the
entropy, we can write $ \eta = \int \eta_x \,d\eta(x) $ and
$h_\eta(f) = \int h_{\eta_x}(f) \,d\eta(x),$ where each $\eta_x$
denotes an ergodic component of $\eta$. Clearly $\eta_x(\La)=1$ for
$\eta$-almost every $x$. Let $\cP$ be a small finite partition of
the space $\cM(\La)$ of invariant probability measures supported in
$\La$ such that
\begin{equation}\label{eq.approx}
\Big| \int \psi_j \,d\xi_1-\int \psi_j \,d\xi_2 \Big| <\beta
\end{equation}
for every $1\le j\le e$ and every pair of probability measures
$\xi_1, \xi_2$ in the same partition element. Set $k=\#\cP$ and
$a_i=\eta(P_i)$ for every element $P_i$ in $\cP$. For every $1\leq i
\leq k$ pick an ergodic measure $\eta_i=\eta_{x_i} \in P_i$
satisfying $h_{\eta_x}(f) \leq h_{\eta_i}(f)+\beta$ for
$\eta$-almost every $\eta_x\in P_i$. 
Part (i) in the lemma is immediate. On the other hand, (ii) follows
because
$$
h_\eta(f)
        = \int h_{\eta_x}(f) \,d\eta(x) \leq \sum_{i=1}^k a_i \,h_{\eta_i}(f)+\beta
        = h_{\hat\eta}(f)+\beta.
$$
Finally, \eqref{eq.approx} implies that
$$
\Big|\int \psi_j \,d\eta - \int \psi_j \, d\hat\eta \Big|
        = \Big| \int \left(\int \psi_j \,d\eta_x\right) \;d\eta(x)
           - \sum_{i=1}^k a_i \int \psi_j \;d\eta_i \Big|
        \leq \sum_{i=1}^k a_i \beta =\beta
$$
for every $j$. This proves (iii) and finishes the proof of the
lemma.
\ep\medbreak

Now we will finish the proof of Theorem~\ref{thm.weakGibbs.strong}.

\proc{Proof}[Proof of Theorem~\ref{thm.weakGibbs.strong}(continuation)]

Take $g\in C(M,\R)$, $c\in \R$, $\psi \in \cF(\La)$ and let $\eta$
be an invariant probability measure such that $\eta(\La)=1$ and
$\int g \,d\eta > c$. Denote by $B_n$ the set of points $x\in M$
such that $S_n g(x)>cn$. Take $\beta>0$ arbitrary small,
$\de_2=\frac{1}{5}(\int g\,d\eta-c)$ and 
the measure $\hat \eta= \sum_{i=1}^k a_i \eta_i$ given by
Lemma~\ref{l.ergodic.approximation} that satisfies
$$
h_{\hat\eta}(f) \geq h_{\eta}(f) - \beta
        \; , \;
\int g \;d\hat\eta \geq \int g \;d\eta -\beta
        \quad\text{and}\quad
\int \psi \;d\hat\eta \leq \int \psi \;d\eta +\beta.
$$
Since $\beta$ is small we can assume $\int g \,d\hat\eta>c+4\de_2$.
Now we claim that
\begin{equation}\tag{$\star\star$}
 \liminf_{n \to \infty}
    \frac{1}{n}\log \nu\left( x\in M :\frac{1}{n} S_n g(x) > c\right)
    \geq h_{\hat\eta}(f) - \int \psi\,d{\hat\eta}-4\beta.
\end{equation}
As before, we may choose $n_0$ sufficiently large and $\de$ small
enough so that, for every $1\leq i\leq k$, the set $D_i$ of points
$x\in \Lambda$ such that
$$
\frac{1}{n} S_n g(x) > \int g \,d\eta_i -\beta,
    \quad
\frac 1 n S_n \psi(x) < \int \psi \,d\eta_i+\beta
    \quad\text{and}\quad
K_n(x,\vep)^{-1} \geq e^{-\beta n}
$$
for every every $n\geq n_0$ and $0<\vep \leq \de$ has
$\eta_i$-measure at least $\frac12$. Hence, given large $n$, small
$\vep>0$ and $1\leq i\leq k$ we proceed as in the proof of
Lemma~\ref{l.lower.ergodic} to obtain a finite set $E_{n}^i
\subset D_i$ 
so that
\begin{enumerate}
\item $E_n^i$ is a maximal $([a_i n],\vep)$-separated set in $D_i$;
\item $\# E_{n}^i \geq e^{ (h_{\eta_i}(f)-\beta) \,[a_i n]}$; and
\item for every $x \in E_n^i$ it holds
        $$
        \frac{1}{[a_i n]} S_{[a_i n]} g(x) > \int g \,d\eta_i-\beta
         \quad\text{and}\quad
        \frac{1}{[a_i n]} S_{[a_i n]} \psi(x) < \int \psi \,d\eta_i+\beta.
        $$
\end{enumerate}
By the specification property, for every sequence
$(x_1,x_2,\dots,x_k)$ with $x_i \in E_n^i$ there exists $x\in M$
that $\vep$-shadows each $x_i$ during $[a_i n]$ iterates with a time
lag of $N(\vep)$ iterates in between. Consequently, if $n$ is large
and $\ti n= \sum_{i} [a_i n] + k N(\vep)$ then $ S_{\ti n}
g(x)>(c+2\de_2) \ti n$. Since the dynamical ball $B(x,\ti n,\vep/8)$
is contained in $B_{\ti n}\cap B(x_1,\ti n, \de_0)$ for every large
$n$, it follows from ~\eqref{eq:abstract.measure.control} that
$$
\nu\Big( B(x,\ti n,\vep) \Big)
    \geq K_{\ti n}(x_1,\vep)^{-1} \;e^{-S_{\ti n}\psi (x)}
    \geq e^{-\beta \ti n} \; e^{-(\int \psi \,d\eta_i+2\beta)\ti n}.
$$
On the other hand, there are at least $\#E_n^1 \times \dots \times
\#E_n^k$ such pairwise disjoint dynamical balls contained in $B_{\ti
n}$. It follows that
$$
\nu(B_{\ti n})
        \geq \sum_{x} \nu(B(x,\ti n,\vep))
        \geq e^{ [ h_{\hat\eta}(f)-\int \psi d\hat\eta-4\beta] \ti n }
$$
for every large $n$, which gives $(\star\star)$. Since
$\beta$ was chosen arbitrarily small 
then
$$
 \liminf_{n \to \infty}
    \frac{1}{n}\log \nu\left( x\in M :\frac{1}{n} S_n g(x) > c\right)
    \geq h_{\eta}(f) - \int \psi\,d{\eta}-6\beta,
$$
which proves the third part in Theorem~\ref{thm.weakGibbs.strong}
and finishes its proof.
\ep\medbreak

%%%%%%%%%%%%%%%%%%%%%%%%%%%%%%%%%%%%%%%%%%%%%%%
\section{Deviation estimates for non-uniformly expanding maps}\label{s.proof2}

In this section we use some of the ideas involved in the proof of
Theorem~\ref{thm.weakGibbs.strong} together with the key notion of
non-uniform specification to prove the large deviation bounds in
Theorem~\ref{thm.weakGibbs.weak}. Through the section, let $M$ be a
compact manifold and $f:M \to M$ be a $C^{1+\al}$ local
diffeomorphism outside a critical/singular region $\cC$ that
satisfies (H). Let $\phi: M\setminus \cC \to \R$ be an H\"older
continuous potential such that (P1)-(P3) hold. Denote by $\nu$ the
corresponding weak Gibbs measure and by $\mu$ the unique equilibrium
state for $f$ with respect to $\phi$.

%%%%%%%%%%%%%%%%%%%%%%%%%%%%%%%%%%%%%%%%%%%%%%
\subsection{Upper bound}

In this subsection we obtain an upper bound for the measure of the
deviation set of non-uniformly expanding maps.

\begin{lemma}\label{l.upperB}
If $g \in C(M,\R)$ and $c \in \R$ then it holds that
\begin{align*}
 \limsup_{n \to \infty} \frac{1}{n}\log \nu & \left(x \in M : \frac{1}{n} S_n g(x) \geq c\right) \\
        & \leq \max\left\{\sup \left\{-P+h_\eta(f) + \int \phi \,d\eta\right\},
        \limsup_{n\to\infty} \frac1n \log \mu(\Gamma_{n})\right\},
\end{align*}
where the supremum is taken over all $f$-invariant measures $\eta$
such that $\int g \,d\eta \geq c$.
\end{lemma}

\proc{Proof}
Let $\be>0$ be given. First we observe that the computations in
Lemma~\ref{le:subexpgrowth} show that $\psi=\phi-P\in \cF(H)$ and
that there exists $C_\beta>0$ such that the set $\De_n$ as in
\eqref{eq.set.Delta.n} is given by
$$
\De_n \supset \{x \in H: n_1(f^{n_i(x)}(x)) \le C_\beta n\}
$$
where $n_i(x) \le n \le n_{i+1}(x)$ are consecutive hyperbolic times
for $x$. In particular, using that $d\mu/d\nu \ge K^{-1}$,
computations analogous to the ones in the proof of
Lemma~\ref{le:subexpgrowth} also give that
\begin{align*} 
\limsup_{n\to\infty} \frac1n \log \nu(\De_n^c)
    & \leq\limsup_{n\to\infty} \frac1n \log \nu\!\left[\bigcup_{1 \leq k \leq n} \{x \in H \!: n_i(x)=k, \, n_1(f^k(x))>C_\beta n\}\right]
   \\ & \leq \limsup_{n\to\infty} \frac1n \log \big(K n \mu (x \in H: n_1(x)>C_\beta n)\big)\\
    & \leq \limsup_{n\to\infty} \frac1n \log \big(K C_\beta^{-1} n \mu (x \in H: n_1(x)>n)\big)\\
    & = \limsup_{n\to\infty} \frac1n \log \mu \left(\Ga_n\right).
\end{align*}
Hence $\de(\vep,\beta)$ does not depend neither on $\vep$ or
$\beta$. Thus, the lemma is an immediate consequence of the
first part in Theorem~\ref{thm.weakGibbs.strong}.
\ep\medbreak

%%%%%%%%%%%%%%%%%%%%%%%%%%%%%%%%%%%%%%%%%%%%%%%%
\subsection{Lower bound estimates}

To obtain lower bound estimates one technical difficulty to overcome
is that no a priori estimates for the measure of dynamical balls hold for specified orbits  
even if the dynamical system satisfies the specification property. 
Here we make use of approximation Lemma~\ref{l.ergodic.approximation} 
and specification properties to prove the following lower bound for the measure of deviation sets.

\begin{proposition}\label{p.lowerB}
Assume that $g \in C(M,\R)$ and $c \in \R$. If either $f$ satisfies the
specification property or $(f,\mu)$ satisfies the non-uniform specification
property then
$$
 \liminf_{n \to \infty} \frac{1}{n}\log \nu\left( x\in M :\frac{1}{n} S_n g(x) > c\right)
    \geq -P+h_\eta(f) + \int \phi\,d\eta,
$$
for every invariant and \emph{expanding} probability measure $\eta$
satisfying $\int g \,d\eta > c$ and $n_1 \in L^1(\eta)$.
\end{proposition}

\proc{Proof}
Note that $\psi=\phi-P \in \cF(H)$ by Lemma~\ref{le:subexpgrowth}.
Set $g \in C(M,\R)$ and $c \in \R$, and
let $B_n$ be the set of points $x\in H$ such that $S_ n g(x)
>cn$. Fix $\beta>0$ arbitrary small and let $\eta$ be an
$f$-invariant and expanding probability measure such that $\int g
\,d\eta > c$ and $n_1 \in L^1(\eta)$. Set also $\de_2=\frac15(\int g
\,d\eta-c)$. Observe that almost every ergodic component $\eta_x$ of
the invariant measure $\eta$ satisfy $n_1 \in L^1(\eta_x)$.
It follows from Lemma~\ref{l.ergodic.approximation} that there are
exists a probability vector $(a_1, \dots, a_k)$ and $f$-invariant
ergodic probability measures $(\eta_i)_{1 \leq j \leq k}$ such that
$\hat\eta=\sum a_j \eta_j$ satisfies
$$
h_{\hat\eta}(f) \geq h_{\eta}(f) - \beta
        \; , \;
\int g \;d\hat\eta \geq \int g \;d\eta -\beta
        \quad\text{and}\quad
\int \psi \;d\hat\eta \leq \int \psi \;d\eta +\beta.
$$
Moreover, it is not hard to check that we can assume $n_1 \in
L^1(\eta_j)$ for every $1 \le j \le k$. So, the Ergodic Theorem and
Remark~\ref{rmk:hyperbolic.times} guarantee that one can pick
$n_0\geq 1$ large and $\de$ small enough such that, for every $1\leq
j\leq k$, the set $D_j$ of points $x\in H$ such that
$$
n_{i+1}(x)-n_i(x) \le \beta n,
    \quad
\frac{1}{n} S_n g(x) > \int g \,d\eta_j -\beta
    \quad\text{and}\quad
\frac 1 n S_n \psi(x) < \int \psi \,d\eta_j+\beta,
$$
for every $n\geq n_0$, has $\eta_j$-measure larger than $\frac12$.
Recall that $n_i(x) \le n <n_{i+1}(x)$ are consecutive hyperbolic
times for $x$. If $0<\vep \ll \de_1$ (as in Lemma~\ref{eq. backward
distances contraction}) is small then $|g(x)-g(y)|<\de_2$ whenever
$|x-y|<\vep$. As in Subsection~\ref{ss.lowerA}, for every large $n$
and small $\vep>0$
there exists a set $E_{n}^j \subset D_j$ such that
\begin{enumerate}
\item $E_n^j$ is a maximal $([a_j n],\vep)$-separated set;
\item $\# E_{n}^j \geq e^{ (h_{\eta_j}(f)-\beta) [a_j n]}$;
\item for every $x \in E_n^j$ it holds
        $$
        \frac{1}{[a_j n]} S_{[a_j n]} g(x) > \int g \,d\eta_j-\beta
         \quad\text{and}\quad
        \frac{1}{[a_j n]} S_{[a_j n]} \psi(x) < \int \psi \,d\eta_j+\beta.
        $$
\end{enumerate}
Condition (1) yield that the dynamical balls $B(x,[a_j n],\vep)$
centered at points in $E_n^j$ are pairwise disjoint. We divide the
remaining of the proof in two cases:

\vspace{.35cm}\noindent \textsf{\bf First case:} $f$ satisfies the
specification property \vspace{.1cm}

Given any sequence $(z_1,z_2,\dots,z_k)$ with $z_j \in E_n^j$ there
exists some point $z\in M$ that $\vep$-shadows each $z_j$ during
$\ell_j:=[a_j n]$ iterates with a time lag of $p_j = N(\vep)$
iterates as in Definition~\ref{d.strong.specification}. Let
$n_j:=n_i(z_j)$ denote the last hyperbolic time for $z_j$ smaller
than $\ell_j$ and write $\ell_j=n_{j}+t_j$ for some $t_j\ge 0$.

\begin{figure}[hbt]
\psfrag{a}{$\sum_{j=1}^{k-1} (\ell_j+p_j)$}
\psfrag{b}{$\sum_{j=1}^{k-1} (\ell_j+p_j)+\ell_k$}
\psfrag{c}{$\sum_{j=1}^{k} (\ell_j+p_j)$}
\psfrag{d}{$n_k$}
\psfrag{3}{$t_k$}
\includegraphics[width=10cm]{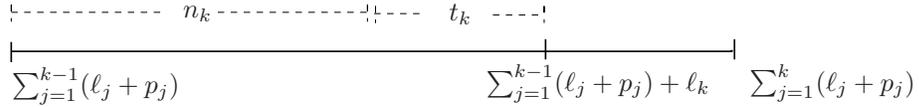}
\caption{Combinatorics of the specified orbit.}
\end{figure}

Therefore, if we set $p_k=0$ and take $\ti n= \sum_{j=1}^k \, \big(
\ell_j + p_j \big)$ one can use that $\max\{k,t_k\} \le \beta \ell_k
\ll \beta (\ti n-t_k)$ to deduce that
\begin{align*}
S_{\ti n} g(z)
        & \geq \sum_{j=1}^{k} S_{\ell_j} \,g(z_j)
              -\de_2 \sum_{j=1}^{k} \ell_j
              - \sup |g| \, k\, N(\vep) \\
        & \geq \sum_{j=1}^k \Big(\int g \,d\eta_j
              - \de_2 -\beta \Big)\, \ell_j
              - \sup |g| \, k\, N(\vep) \\
        & >   (c+3\de_2) n
              - \sup |g| \, k\, N(\vep) \\
        & \geq (c+3\de_2) \ti n
              - 2\sup |g| \, k\, N(\vep)\\
        & \geq (c+2\de_2) \ti n
\end{align*}
provided that $n$ is large enough. Hence $B (z, \ti n, \vep) \subset
B_{\ti n}$ and there are at least $e^{(h_{\hat \eta}(f)-2\beta)\ti
n}$ such distinct dynamical balls. Since each of the points $x_i$
were chosen in a full measure set then the weak Gibbs property yields
an estimate for the measure of the corresponding dynamical balls.
However, no a priori estimates on the measure of the specified orbit $z$ is
guaranteed. We claim that
\begin{align}\label{eq:measure.concatenation}
\nu(B(z,\ti n,\vep))
          & \geq e^{-2k \sup|\psi|\, \beta \ti n}\; \Bigg(\prod_{j=1}^{k} e^{-S_{\ell_i}\psi(z_j)} \Bigg)
          \nonumber
        \\& \geq e^{-2k \sup|\psi|\, \beta \ti n}\, e^{-2\beta \ti n}\, e^{-\ti n \,\int \psi \,d\eta}
\end{align}
for every large $n$.
 
\proc{Proof}[Proof of the claim:]
Let $L$ and $\ga$ be given by condition (C).
For notational simplicity set $\ti z_{k+1}\!=\!f^{\sum_{j=1}^k (\ell_j+p_j)}(z)\in
B(z_{k+1},\ell_{k+1},\vep)$. Since $B(z, \ti n, \vep)$ contains  
$B(z,\ti n+n_{k+1}, \vep)$, where $\ell_k+n_{k+1}$ denotes the first
hyperbolic time for $z_k$ larger than $\ell_k$, first we show that
\begin{equation}\label{eq:big.ball}
f^{\ti n+n_{k+1}}(B(z,\ti n+n_{k+1}, \vep))
    = B(f^{\ti n+n_{k+1}}(z),\vep).
\end{equation}
Since $\ell_k+n_{k+1}$ is a hyperbolic time for $z_{k}$ and $\vep
\ll \de_1$ then there exists backward distance contraction and
$f^{\ell_k+n_{k+1}} (B(\ti z_{k},\ell_k+n_{k+1},2\vep))
        = B(f^{\ell_k+n_{k+1}}(\ti z_{k}),2\vep)$.
In fact, using $d(f^{n_{k}}(\ti z_{k}), f^{n_{k}}(z_{k}))<\vep$, that $\beta$ is fixed arbitrary small and 
$\ell_k-n_k<\beta n$ (recall the definition of the set $D_k$) then
$$
 \diam(B(\ti z_{k},\ell_k+n_{k+1},\vep))
        \le \diam(B(\ti z_{k},n_{k},\vep))
        \le \vep \si^{-\frac12 n_{k}}
        \le \vep \si^{-\frac12 [a_k-\beta]n} 
        \ll \vep
$$
provided that $n$ is large enough. Again, if $n$ is large, using property (C) on the diameter
of preimages this yields that there exists $\ti L>0$ so that 
\begin{align*}
 \diam(f^{-p_{k-1}-t_{k-1}}(B(\ti z_{k},\ell_k+n_{k+1},\vep)))
        & \le \ti L [\vep \si^{-\frac12 [a_k-\beta]n}]^{\ga^{p_{k-1}+t_{k-1}}}\\
        & \le \ti L [\vep \si^{-\frac12 [a_k-\beta]n}]^{\ga^{N(\vep)+\beta n}}
        \ll \vep
\end{align*}
and, consequently,
$
f^{-p_{k-1}-t_{k-1}}(B(\ti z_{k},\ell_k+n_{k+1},\vep))
        \subset B(f^{n_{k-1}}(\ti z_{k-1}),2\vep).
$
Recall that $\ti z_{k}=f^{\ell_{k-1}+p_{k-1}}(\ti z_{k-1})$, that
$n_{k-1}$ is a hyperbolic time for $z_{k-1}$ and there exists backward distance
contraction in the dynamical ball of radius $\de_1 \gg 2\vep$. Hence
$
f^{n_{k-1}} (B(\ti z_{k-1},n_{k-1},2\vep))
        = B(f^{n_{k-1}}(\ti z_{k-1}),2\vep)
$ 
and
\begin{align*}
 B(\ti z_{k-1}, \ell_{k-1} & +p_{k-1}+ \ell_k+n_{k+1},\vep)\\
        &= B(\ti z_{k-1},n_{k-1},\vep) \cap f^{-\ell_{k-1}-p_{k-1}}(B(\ti z_{k},\ell_k+n_{k+1},\vep))\\
        &= B(\ti z_{k-1},n_{k-1},\vep) \cap f^{-n_{k-1}}[f^{-p_{k-1}-t_{k-1}}(B(\ti z_{k},\ell_k+n_{k+1},\vep))]\\
        &\subset  B(\ti z_{k-1},n_{k-1},\vep) \cap f^{-n_{k-1}}[B(f^{n_{k-1}}(\ti z_{k-1}),2\vep)].
\end{align*}

\begin{figure}[hbt]
\psfrag{a}{$z_j$}
\psfrag{b}{$f^{n_j}(z_j)$}
\psfrag{c}{$B(f^{n_j}(z_j),2\vep)$}
\psfrag{e}{$\ti z_{j+1}$}
\psfrag{d}{$B(f^{n_{j+1}}(\ti z_{j+1}),2\vep)$}
\includegraphics[width=10cm]{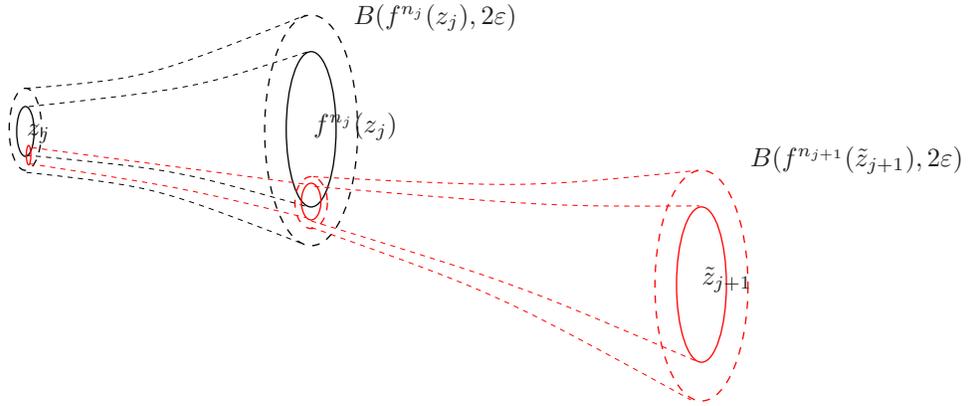}
\caption{Concatenation of dynamical balls. Color online}
\end{figure}

Consequently, the dynamical ball $B(\ti z_{k-1},
\ell_{k-1}+p_{k-1}+\ell_k+n_{k+1},\vep)$ is mapped diffeomorphically
by $f^{\ell_{k-1}+p_{k-1}+\ell_k+n_{k+1}}$ onto the ball centered at
$f^{\ell_k+n_{k+1}}(\ti z_{k})$ 
with radius $\vep$. 
Using the same argument as above recursively we obtain \eqref{eq:big.ball} as desired.

It remains to compute the measure of $B(z,\ti n + n_{k+1},\vep)$.
Using \eqref{eq:big.ball}, the fact that $\nu$ is a conformal
measure with Jacobian $J_\nu f= e^{-\psi}$ and the bounded
distortion property at hyperbolic times (see Corollary~\ref{c.bounded.distortion})
it follows that
\begin{align*}
\nu(B(f^{\ti n + n_{k+1}}(z),& \vep))
          = \int_{B(z,\ti n + n_{k+1},\vep)} e^{S_{\ti n + n_{k+1}}\psi(y)} \; d\nu(y)\\
         & \le K_0^k  e^{\sup|\psi|\,[n_{k+1}+\sum_{j} (t_j + p_j)]} \,\Bigg(\prod_{j=1}^{k} e^{S_{\ell_i}\psi(z_j)} \Bigg)
                    \, \nu(B(z,\ti n + n_{k+1},\vep)).
\end{align*}
Since $\nu$ is an open measure then every ball of radius $\vep$ has
measure at least $C_\vep>0$. Finally, using that $n_k \le  \ell_k \le (\ell_k+n_{k+1})$ are 
consecutive hyperbolic times for $z_k$ then $n_{k+1} \le (\ell_k+n_{k+1})-n_k 
\le \beta \ell_k=\beta[a_k n] \ll \beta \ti n$ and so
\begin{align*}
\nu(B(z,\ti n,\vep)) & \geq \nu(B(z,\ti n + n_{k+1},\vep)) \\
        & \geq C_\vep K_0^{-k}\; \Bigg(\prod_{j=1}^{k} e^{-S_{\ell_i}\psi(z_j)} \Bigg)\; e^{-\sup|\psi|\,[k N(\vep) + (k+1) \beta \ti n]}\\
        & \geq  e^{-2k \sup|\psi|\, \beta \ti n}\; \Bigg(\prod_{j=1}^{k} e^{-S_{\ell_i}\psi(z_j)} \Bigg)\\
        & \geq e^{-2k \sup|\psi|\, \beta \ti n}\; \exp\Big( \sum_{j=1}^{k} (-\int \psi \,d\eta_j-\beta) \ell_j\Big)\\
        & \geq e^{-2k \sup|\psi|\, \beta \ti n}\, e^{-2\beta \ti n}\, e^{-\ti n \,\int \psi \,d\eta}
\end{align*}
for every large $n$, which proves our claim.
\ep\medbreak

We are now in a position to finish the proof of the first case of
the proposition. Indeed, note that we obtain as a direct consequence
of equation~\eqref{eq:measure.concatenation} that
 $
 \log \nu(B_{ \ti n})
    \geq \Big( h_{\eta}(f)-\int \psi d\eta -5\beta -2\beta k \sup|\psi|\Big) \ti n
 $
for every large $n$. Since $\beta$ was arbitrary this shows that
$$
\liminf_{n \to \infty}
        \frac{1}{n}\log \nu\left( x \in M : \frac{1}{n} S_n g(x) > c  \right)
        \geq h_{\eta}(f) - \int \psi \,d\eta.
$$

\vspace{.55cm}\noindent \textsf{\bf Second case:} $(f,\mu)$
satisfies the non-uniform specification property \vspace{.1cm}

In the case that $(f,\mu)$ satisfies the non-uniform
specification property the computations are similar to the previous
ones with the difference that the time lags given by non-uniform specification may be unbounded.
Take $n_0$ large  and $\de$ small so that
\begin{equation*}
p(x,n,\vep) \leq \beta n 
\end{equation*}
for every $x \in D_i$, $0<\vep \leq \de$, $1 \le
i \le k$ and $n \geq n_0$. Through the remaining of the proof 
 set also $p_j:= \max_{x \in E_n^j} p(x,n,\vep)$.

For every sequence $(z_1,z_2,\dots,z_k)$ with $z_j \in E_n^j$ there
exists some point $z\in M$ that $\vep$-shadows, in the non-uniform
metric, each $z_i$ during $\ell_j:=[a_j n]$ iterates with a time lag
of $p_j$ iterates.  Moreover, if $\ti n= \sum_{j=1}^k
\, \big( \ell_j + p_j \big)$, the set of points $z$ obtained as
above are $(\ti n,\vep)$ separated and there are at least
$e^{(h_{\hat \eta}(f)-2\beta)\ti n}$ such points.
Since $\beta>0$ is small, observe that
\begin{equation}
S_{\ti n} g(z)
        \geq \sum_{j=1}^k S_{\ell_j} \,g(z_j)
              -\de_2 \sum_{j=1}^k \ell_j
              - \sup |g| \sum_{j=1}^k p_j, 
\end{equation}
which is bounded from below by
\begin{equation*}
    \sum_{j=1}^k \Big(\int g \,d\eta_j-2\beta - \de_2 \Big) \, \ell_j
        > \Big(\int g \,d\eta-2\de_2\Big) n
        > (c+ 2\de_2 ) n + (c+\de_2) \frac{\beta n}{\sup |g|}
        > (c+\de_2 ) \ti n
\end{equation*}
for every large $n$. It follows from our choice of $\vep$ that $B
(z, \ti n, \vep) \subset B_{\ti n}$. Given $z$ as above, the
computations involved in the proof of
\eqref{eq:measure.concatenation} give that
\begin{align*}
\nu(B(z,\ti n + n_{k+1},\vep))
        & \geq C(\vep)^{-1} K_0^{-k} e^{-\sup|\psi|\,[n_{k+1}+\sum_{j} (t_j + p_j)]} \,\Bigg(\prod_{j=1}^{k} e^{-S_{\ell_i}\psi(z_j)} \Bigg)\\
        & \geq e^{-3(k+1)\sup|\psi|\,\beta \ti n} \,\Bigg(\prod_{j=1}^{k} e^{-S_{\ell_i}\psi(z_j)} \Bigg)
\end{align*}
and, consequently,
$
\nu(B(z,\ti n,\vep))
         \geq e^{-3(k+1) \sup|\psi|\, \beta \ti n}\, e^{-2\beta \ti n}\, e^{-\ti n \,\int \psi \,d\eta}
$ for every large integer $n$. Henceforth,
 $
 \log \nu(B_{ \ti n})
    \geq \Big( h_{\eta}(f)-\int \psi d\eta -5\beta -3(k+1)\beta \,\sup|\psi| \Big) \ti n
 $
for large $n$. Since both $0<\vep<\de$ and $\beta>0$ were chosen
arbitrarily small and
$
\lim_{\vep \to 0} \limsup_{n \to \infty} \frac{p(x,n,\vep)}{n}=0
$
for almost every $x$ one obtains
$$
\liminf_{n \to \infty}
        \frac{1}{n}\log \nu\left( x \in M : \frac{1}{n} S_n g(x) > c  \right)
        \geq h_{\eta}(f) - \int \psi \,d\eta.
$$
The proof of the proposition is now complete.
\ep\medbreak

%%%%%%%%%%%%%%%%%%%%%%%%%%%%%%%%%%%%%%%%%%%%%%%
\section{Some applications}\label{s.applications}

%%%%%%%%%%%%%%%%%%%%%%%%%%%%%%%%%%%%%%%%%%%%%%%
\subsection{One-dimensional examples}

Large deviations estimates for one-dimensional non-uniformly expanding
maps were obtained only by Keller and Nowicki \cite{KN92} 
for quadratic maps satisfying the Collet-Eckmann condition and by
Ara\'ujo and Pac\'ifico \cite{AP06} to non-uniformly expanding
quadratic maps. The first authors proved a large deviations
principle for observables of bounded variation and the second
authors obtained upper bounds for the measure of the deviation sets
of any continuous observable.
Using that every topologically mixing and continuous interval map satisfies specification
(see ~\cite{Bl83}) we will now discuss applications of our
results to some important classes of examples.

\begin{example}\emph{(Non-uniformly expanding quadratic maps)}\\
We consider the class of quadratic maps $f_a$ on the real line 
given by $$f_a(x)=1-ax^2.$$ In \cite{BC85}, Benedicks
and Carleson proved the existence of a positive Lebesgue measure set
of parameters $\Omega \in [0,2]$ such that for every $a\in \Omega$
the quadratic map $f_a$ has positive Lyapunov exponent and an unique
absolutely continuous invariant probability measure $\mu_a$
supported on $[f^2(0),f(0)]$.  In fact, these maps are topologically mixing on $[f^2(0),f(0)]$
and $d\mu_a/d\Leb \in L^p$ for every $p<2$. 
It follows from the previous discussion that each $f_a$
satisfies the specification property.
Moreover, the same argument used in \cite{AP09} to deal with
infinitely many critical points is enough to guarantee that $\Leb(\Gamma_n)$
decays exponentially fast (cf. \cite[Section~2.1]{AP06}).

Now we notice that all invariant measures are expanding. Indeed,  on the
one hand \cite[Proposition~3.1]{BK98} establishes for $S$-unimodal maps and any invariant measure 
$
\lambda(\mu) \geq \lambda_{\text{per}},
$
where $\lambda(\mu)=\int \log|f'| \, d\mu$  is the integrated Lyapunov exponent of $\mu$
and $\lambda_{\text{per}}$ is the infimum of Lyapunov exponents among periodic orbits. On the
other hand, it follows from \cite{NS98} that the Collet-Eckmann is
equivalent to $\lambda_{\text{per}}>0$.

Since there exist many invariant probability measures with integrable first
hyperbolic time map we proceed to show that the measure of the deviation
sets is exponential. Using that $d\mu_a/d \Leb \in L^p$ for any $p\in(1,2)$ and that $\Leb(\Gamma_n)$ decreases exponentially fast then, if  $q> 1$ satisfies $\frac1p+\frac1q=1$, by  H\"older's inequality
$$
\mu_a(\Gamma_n)
	=\int 1_{\Gamma_n} d\mu_a
	 = \int 1_{\Gamma_n} \frac{d\mu_a }{d\nu} d\nu
	 \le \left\|\frac{d\mu_a}{d\nu}\right\|_p \, Leb(\Gamma_n)^q
$$
also decreases exponentially fast and $n_1\in L^1(\mu_a)$. Since $\mu_a$ is an equilibrium state for $\phi_a=-\log |f_a'|$ then it follows from Ruelle-Pesin's formulas that $P=0$. Moreover, $\nu=\Leb$ is an expanding conformal measure and so
$$
 \limsup_{n\to \infty} \frac1n\log Leb \left(x \in M : \left|\frac1n S_n g(x)-\!\int g\,d\mu_a\right| \ge c\right)
        \leq -\alpha
$$
where
$$
\alpha=\min\left\{ 
		-\lim_{n\to\infty} \frac1n \log \mu_a(\Gamma_n), 
		\sup\{ -h_\eta(f)+\int \log |f_a'| \,d\eta \}
	\right\} >0,
$$
and the supremum in the right hand term is over all invariant measures $\eta$ such that
$|\int g d\eta - \int g d\mu_a | \ge c$. Analogously,
$$
 \liminf_{n\to \infty} \frac1n\log Leb \left(x \in M : \left|\frac1n S_n g(x)-\!\int g\,d\mu_a\right| > c\right)
        \geq -\beta
$$
where
$
0<\beta= \sup\{ -h_\eta(f)+\int \log |f_a'| \,d\eta\}
$
and the supremum is taken over all invariant measures $\eta$ such that $n_1\in L^1(\eta)$ and $|\int g d\eta - \int g d\mu_a | > c$.
\end{example}

In the following examples we obtain some large deviation estimates for maps of the interval with intermitency
behaviour with respect to some equilibrium states. In particular we consider the case of the physical  
and the maximal entropy measure. Moreover, we discuss the presence of the condition on the decay of the first 
hyperbolic time map in the large deviations upper bound.  

\begin{example}\label{ex.maneville}\emph{(Intermittency phenomena)}\\
Given $\al \in (0,1)$, let $f:[0,1]\to [0,1]$ be the $C^{1+\alpha}$
transformation of the interval given by
\begin{equation*}\label{eq. Manneville-Pomeau}
f_\al(x)= \left\{
\begin{array}{cl}
x(1+2^{\alpha} x^{\alpha}) & \mbox{if}\; 0 \leq x \leq \frac{1}{2}  \\
2x-1 & \mbox{if}\; \frac{1}{2} < x \leq 1.
\end{array}
\right.
\end{equation*}
known as the Maneville-Pomeau map. This transformation has $0$ as an
indifferent fixed point (that is $Df(0)=1$) and expansion everywhere
else. The map presents an intermittency phenomenon. We provide bounds for the 
measure of deviation sets in the context of the SRB measure and the case of the 
maximal entropy measure.
\vspace{.1cm}

\noindent \emph{(a) SRB measure}

It is known that $f$ has a finite absolutely continuous invariant probability
measure $\mu$ with polynomial decay of correlations of order
$\mathcal O(n^{\frac1\alpha-1})$. In fact that is also the decay of
the tail of the first hyperbolic time with respect to $m=\Leb$.

By Ruelle-Pesin's formula, $\mu$ is an equilibrium state for $f$
with respect to the potential $\phi=-\log|Df|$ with pressure
$P:=P(\phi)=0$. In fact, $\mu$ and $\de_0$ are the unique ergodic
equilibrium states for $\phi$, and so any other equilibrium state is of the
form $t\mu + (1-t) \de_0$ for some $t\in(0,1)$. In addition, it is proved
in~\cite{Hu04} that $d\mu/dm \approx x^{-\alpha}$. However, $n_1
\not\in L^1(\mu)$. Roughly, partitioning the unit interval according
to the sequence $(\frac1n)_n$ it follows that
\begin{align*}
\int n_1 d\mu
    & \geq \sum_{n \geq 1} n_1(\frac1{n+1})\mu([\frac1{n+1},\frac1n])
    \approx \sum_{n \geq 1} n_1(\frac1{n+1})m([\frac1{n+1},\frac1n])(\frac1n)^{-\alpha}
    \\ &= \sum_{n \geq 1} n^\alpha \;n_1(\frac1{n+1})m([\frac1{n+1},\frac1n]),
\end{align*}
which is infinite because $n_1\ge 1$. 
In consequence the Lebesgue measure of deviation sets decrease polynomially and
$
\limsup \frac1n\log\mu(\Ga_n)=0.
$
Since $f$ admits a finite and generating Markov partition then it satisfies the specification property.
Therefore, it follows from Theorem~\ref{thm.weakGibbs.strong} that for every \emph{continuous}
observable $g$ 
\begin{equation}\label{eq:Maneville.Pomeau}
 \liminf_{n\to \infty} \frac1n\log \Leb \left[x \in M : \left|\frac1n S_n g(x)-\!\int g\,d\mu\right|>c\right]
        \geq \sup_{\eta} \left\{h_\eta(f) - \!\int \log|Df| \,d\eta\right\}
\end{equation}
where $\eta$ denotes an invariant measure so that $\eta(H)=1$,
$|\int g \,d\eta-\int g\,d\mu|>c$ and $n_1 \in L^1(\eta)$, and
\begin{align}\label{eq:Maneville.Pomeau2}
 \limsup_{n\to \infty} \frac1n & \log \Leb \left(x \in M : \left|\frac1n S_n g(x)-\!\int g\,d\mu_\phi\right| 
 	\ge c\right) \nonumber \\
        & \leq \max
        \left\{
            \sup \left\{-P+h_\eta(f) + \int \phi \,d\eta\right\},
            \limsup_{n\to\infty} \frac1n \log \mu_\phi(\Gamma_n)
        \right\}
\end{align}
where the supremum is taken over all invariant probability measures $\eta$ such that 
$|\int g \,d\eta -\int g d\mu_\phi|\geq c$.
In the case that $g(0)=\int g \,d\mu$ then any invariant probability measure $\eta$ considered in the 
right hand side of \eqref{eq:Maneville.Pomeau} and \eqref{eq:Maneville.Pomeau2} is far from the 
convex hull generated by the equilibrium states $\Leb$ and $\de_0$. Hence it holds that the supremum over
all invariant probability measures $\eta$ such that  $|\int g \,d\eta -\int g d\mu_\phi|\geq c$ satisfies
$$
\sup \left\{-P+h_\eta(f) + \int \phi \,d\eta\right\} < 0
	=  \limsup_{n\to\infty} \frac1n \log \mu_\phi(\Gamma_n)
$$
and that the measure of deviation sets decrease at most exponentially fast. Some results in \cite{AFLV10}
relate decay of correlations, the decay of the tail of inducing maps and the rate of decay of the deviations with 
respect to the invariant probability measure. The tail of the first hyperbolic time in \eqref{eq:Maneville.Pomeau2} 
gives an indication that this relation can be expected to hold also in the case of deviations with respect to the
Lebesgue measure. 
Let us also point out that the results obtained by Chung~\cite{Chu11} yield a large deviations principle, where 
the rate function is not strictly concave due to the non-uniqueness of equilibrium states. In particular, for an open
and dense set of observables (namely those satisfying $g(0)\neq\int g \,d\mu$) it follows that deviations are sub-exponential. In fact, we note that that polynomial upper and lower bounds for \emph{H\"older continuous} observables have been established in~\cite{MN08,Mel09, PS09}.
\vspace{.1cm}

\noindent \emph{(b) Equilibrium states for potentials with small variation}

It was obtained in \cite{VV10} that $f$ admits a unique equilibrium state $\mu_{\phi}$ with respect to
any H\"older continuous potential $\phi$ such that $\sup\phi-\inf\phi<\log 2$. Moreover, $\mu_\phi$ is absolutely continuous with respect to a weak Gibbs conformal measure $\nu_\phi$, is expanding and $\mu_\phi(\Ga_n)$ 
decays exponentially fast. Hence it follows from Corollary~\ref{cor:conseq} that for every continuous observable $g$
it holds that
\begin{align*}
 \limsup_{n\to \infty} \frac1n & \log \nu_\phi \left(x \in M : \left|\frac1n S_n g(x)-\!\int g\,d\mu_\phi\right| 
 	\ge c\right)\\
        & \leq \max
        \left\{
            \sup \left\{-P+h_\eta(f) + \int \phi \,d\eta\right\},
            \limsup_{n\to\infty} \frac1n \log \mu_\phi(\Gamma_n)
        \right\}
\end{align*}
where the supremum is taken over all invariant probability measures $\eta$ such that 
$|\int g \,d\eta -\int g d\mu_\phi|\geq c$, and also that 
\begin{align*}
\liminf_{n\to \infty} \frac1n & \log \nu_\phi \left(x \in M : \left|\frac1n S_n g(x)-\!\int g\,d\mu_\phi\right|>c\right)\\
        &\geq \sup \left\{-P+h_\eta(f) + \int \phi\,d\eta\right\}
\end{align*}
where the supremum is taken over all invariant
probability measures $\eta$ such that
$\eta(H)=1$, $|\int g \,d\eta -\int g d\mu_\phi|>  c$ and
$n_1 \in L^1(\eta)$. Since the equilibrium state is unique the right hand side of both expressions above
is strictly negative, which yields that the measure of deviation sets decrease exponentially fast.
We also remark that if $\eta$ is an $f$-invariant probability measure with $\eta\neq \de_0$ 
then it follows from Birkhoff's ergodic theorem that its Lyapunov exponent is $\int \log |f'| d\eta>0$. Therefore
for an open and dense class of continuous observables $g$ (namely those that  satisfy $g(0)\neq \int g\,d\mu_\phi$) 
if $c$ is small enough then
\begin{align*}
 \limsup_{n\to \infty} \frac1n  & \log  \nu_\phi \left(x \in M :  \left|\frac1n S_n g(x)-\!\int g\,d\mu_\phi\right| 
 	\ge c\right)\\
          & \leq \sup \left\{-P+h_\eta(f) + \int \phi \,d\eta\right\}
\end{align*}
where the supremum is taken over all invariant expanding probability measures $\eta$ such that
$|\int g \,d\eta -\int g d\mu_\phi| \ge  c$. This indicates that it may be possible to establish a large deviations 
principle for this non-uniformly expanding dynamics using expanding measures.  We refer the reader to
Section~\ref{s.some.results} for further discussion.
\end{example}

%%%%%%%%%%%%%%%%%%%%%%%%%%%%%%%%%%%%%%%%%%%%%
\subsection{Higher dimensional examples}\label{s.higher.dimensions}

The next class of examples are multidimensional local
diffeomorphisms obtained by local bifurcation of expanding maps and
were introduced in \cite{ABV00}. Although the original expanding
maps satisfy the specification property we point out that the same
should not hold for the perturbations.

\begin{example}
Let $f_0$ be an expanding map in $\mathbb T^n$ and take a periodic
point $p$ for $f_0$. Let $f$ be a $C^1$-local diffeomorphism
obtained from $f_0$ by a bifurcation in a small neighborhood $U$ of
$p$ in such a way that:
\begin{enumerate}
\item every point $x \in \mathbb T^n$ has some preimage outside $U$;
\item $\|Df(x)^{-1}\| \leq \si^{-1}$ for every $x \in \mathbb T^n\backslash
U$, and $\|Df(x)^{-1}\| \leq L$ for every $x \in \mathbb T^n$ where $\si>1$ is
large enough or $L>0$ is sufficiently close to $1$;
\item $f$ is topologically exact: for every open set $U$ there is $N\geq 1$ for which $f^N(U)=\mathbb T^n$
\end{enumerate}
It follows from  \cite{VV10} that $f$ has a unique
(ergodic) equilibrium state $\mu$ for the H\"older continuous
potential $\phi=-\log |\det Df|$, it is absolutely continuous with
respect to the conformal measure $\nu=\Leb$ with density bounded
away from zero and infinity, and it is expanding. 
We note also that the equilibrium state $\mu$ also satisfies
the non-uniform specification property.

\begin{lemma}\label{le:local}
$(f,\mu)$ satisfies the non-uniform specification property.
\end{lemma}

\proc{Proof}
First we note that
since $M$ is compact and $f$ is topologically exact then for every 
$\vep>0$ there exists $N_\vep\ge 1$ such that $f^{N_\vep}(B)=\mathbb T^n$ 
for every ball $B$ of radius $\vep$. Indeed, for every $x$ let $N(x,\vep)\ge 1$
be the minimum integer such that  $f^{N(x,\vep)}(B(x,\vep/3))=\mathbb T^n$. By compacteness
the open cover $(B(x,\vep/3))_{x\in \mathbb T^n}$ admits a finite covering $(B(x_i,\vep/3))_{i=1..n}$.
Hence, if $N_\vep=\max\{N(x_i,\vep): i=1..n\}$ then 
any ball $B$ of radius $\vep$ contains a ball $B(x_j,\vep/3)$, for some $j$, and so
$f^{N_\vep}(B)=\mathbb T^n$.

It follows from~\cite{VV10} that the equilibrium state
$\mu$ is absolutely continuous with respect to a conformal measure
$\nu$ with density bounded away from zero and infinity and $n_1 \in
L^1(\mu)$. Moreover, 
the sequence $n_k(\cdot)$ of hyperbolic times is non-lacunar, that is
$\frac{n_{k+1}-n_k}{n_k} \to 0$ at almost every $x$. Therefore,
if $0<\vep<\de$, $n$ is large and $n_k(x)<n<n_{k+1}(x)$ are consecutive hyperbolic times
then clearly $B(x, n_{k+1},\vep) \subset B(x,n,\vep)$ and 
$$
f^{n_{k+1}+N_\vep}(B(x,n_{k+1},\vep))
	= f^{N_\vep}(  B(f^{n_{k+1}}(x),\vep)  ) 
	= \mathbb T^n.
$$
Thus for any given  $y\in \mathbb T^n$ and proximity $\zeta>0$ there exists 
$z\in B(x,n,\vep)$ so that 
$f^{N_\vep+n_{k+1}(x)-n}(f^n(z))=f^{N_\vep+n_{k+1}(x)}(z) \in B(y,\zeta)$.  
Take  $p(x,n,\vep)=N_\vep+ n_{k+1}(x)-n$. Then for any
$x_1, \dots, x_m$ in a full $\mu$-measure set, any positive integers $k_1, \dots, k_m$ and
$p_i\ge p(x_i,n_i,\vep)$ there exists $z\in \mathbb T^n$ such that 
$
z \in B(x_1,n_1,\vep)
$
and
$
f^{n_1+p_1+\dots +n_{i-1}+p_{i-1}}(z) \in B(x_i,n_i,\vep)
$
for every $2\leq i\leq k$. To obtain the non-uniform specification property just note that
$$
\lim_{\vep\to 0} \limsup_{n\to \infty} \frac{p(x,n,\vep)}{n}
	    \leq \lim_{\vep\to 0}\limsup_{k\to \infty} \frac{N_\vep+ n_{k+1}(x)-n_k(x)}{n_{k}(x)}=0.
$$
\ep\medbreak

Using the non-uniform specification property we obtain from Theorem~\ref{thm.weakGibbs.weak}
that
\begin{align*}\label{eq:local.diffeo.I}
 \limsup_{n\to \infty} \frac1n & \log \Leb \left(x \in M : \left|\frac1n S_n g(x)-\!\int g\,d\mu\right| \geq c\right)\\
        & \leq \max
        \left\{
            \sup \left\{-P+h_\eta(f) + \int \phi \,d\eta\right\},
            \limsup_{n\to\infty} \frac1n \log \mu(\Gamma_n)
        \right\},
\end{align*}
where the supremum taken over all invariant probability measures
$\eta$ satisfying $|\int g \,d\eta -\int g \,d\mu|>c$, and also
\begin{align*}
 \liminf_{n\to \infty} \frac1n & \log \Leb \left(x \in M : \left|\frac1n S_n g(x)-\!\int g\,d\mu\right|>c\right)\\
        & \geq \sup \left\{-P+h_\eta(f) + \int \phi\,d\eta\right\},
\end{align*}
where the supremum taken over expanding $f$-invariant probability
measures $\eta$ such that $n_1 \in L^1(\eta)$ and $|\int g \,d\eta
-\int g \,d\mu|>c$. Note that both rates are exponential.
\end{example}

We also prove that a robust class of multidimensional non-uniformly
expanding maps with singularities also satisfy this weak form of specification. 

\begin{example}\emph{(Viana maps)}\\ 
In \cite{Vi97}, the author introduced a robust class of
multidimensional non-uniformly hyperbolic maps with singularities
commonly known as Viana maps. More precisely, these are obtained as
$C^3$ small perturbations of the skew product $\phi_\al$ of the
cylinder $S^1 \times I$ given by
\begin{equation*}\label{eq:Viana.maps}
\phi_\al(\theta,x)
    =(d \theta (\hspace{-.3cm}\mod 1)\;,\; 1-ax^2+\alpha \cos(2\pi \theta)),
\end{equation*}
where $d\geq 16$ is an integer, $a$ is a Misiurewicz parameter for
the quadratic family, and $\alpha$ is small. These maps admit a
unique SRB measure $\mu$ (it is absolutely continuous with respect
to $m=\Leb$, has only positive exponents and $d\mu/dm \in L^p(m)$
where $p=d/(d-1)$) and are strong topologically mixing on the
attractor $\La=\cap_{n \geq 0} \phi_\al^n(S^1\times I)$: for every open set
$A$ there exists a positive integer $n=n(A)$ such that $\phi_\al^n(A)=\La$.
See~\cite{Vi97,Al00,AV02} for more details.

We show that $(f,\mu)$ satisfies the non-uniform specification using
that the sequence of hyperbolic times is non-lacunar, that is,
$\frac{n_{k+1}-n_k}{n_k} \to 0$ at almost every $x$ and that the
image of hyperbolic balls grow to $\La$ after finitely many
iterates.

First we observe that $n_1$ is integrable with respect to the SRB
measure $\mu$. In fact, since the tail of the first hyperbolic time
map decays subexponentially fast with respect to $m$ in particular
one has $n_1 \in L^q(m)$ for every $q \ge 1$.
Using once more Cauchy-Schwartz inequality and 
that $d\mu/dm \in L^p(m)$ it  follows that $n_1 \in L^1(\mu)$ and, consequently, 
the sequence $n_j(\cdot)$ of hyperbolic times is non-lacunar (see e.g. \cite[Corollary~3.8]{VV10}).

The arguments used in the proof of Theorem~C in \cite{AV02}
give that for every image of a rectangle at a hyperbolic time grow
to $\La$ after a finite number $\ell$ of iterates. Since $\Leb$ almost every point
has infinitely many hyperbolic times then $f$ is topologically exact.
Therefore, the argument that $(f,\mu)$ satisfies the non-uniform specification property
goes along the same lines used in proof of Lemma~\ref{le:local}.

For completeness, let us mention that during the refereeing process it has been announced
in \cite{AFLV10} that Viana maps have at least stretched exponential large deviations with 
respect to the invariant SRB measure $\mu$. 
\end{example}

%%%%%%%%%%%%%%%%%%%%%%%%%%%%%%%%%
\section{Recent developments and some future perspectives}\label{s.some.results}

Both weak specification properties and the theory of large deviations
have been target of recent intense study. In this section we discuss 
some recent results and establish a connection with some future perspectives.    

%%%%%%%%%%%%%%%%%%%%%%%
\subsubsection*{Weak specification properties}

The important notion of (strong) specification introduced by Bowen in~\cite{Bo71} allowed to deduce
that uniformly hyperbolic maps are rich from the ergodic theory viewpoint. In fact, this property play a key role 
in the proof that these maps have a unique equilibrium state for every H\"older continuous potential, 
that the topological entropy coincides with the exponential growth rate of the set $Per_n(f)$ of 
periodic points of period $n$. Moreover, it  has important connections with the study of Poincar\'e recurrence 
and large deviations.   

Hence, it is important to understand which topological or measure-theoretical weaker forms of specification 
do non-uniformly hyperbolic dynamical systems satisfy, which ergodic properties are obtained as 
consequence and the relation between the topological and measure-theoretical notions. 
Firstly let us mention that Oliveira~\cite{Ol11} proved for $C^1$-endomorphisms that every ergodic and 
invariant probability measure  with only  positive Lyapunov exponents satisfies the nonuniform specification property
introduced in \cite{STV03} and deduced interesting results concerning Poincar\'e recurrence. Very recently,  Oliveira and Tian~\cite{OT11} announced that every ergodic hyperbolic measure preserved by a $C^{1+\al}$-diffeomorphism satisfy both the nonuniform specification property of \cite{STV03} as the one introduced here. 
In consequence, since Dirac measures satisfy the nonuniform specification property trivially, all invariant 
probability measures for the Maneville-Pommeau transformation in Example~\ref{ex.maneville} do satisfy
the non-uniform specification property. In fact this is already a consequence from the fact that the 
Maneville-Pommeau transformation satisfies the (strong) specification property.  So, having in mind the results
obtained in \cite{AAS03,Cao03} it would be very interesting to answer the following question:

\vspace{.2cm}
{\bf Question 1:} \textit{ Let $f: M\to M$ be a $C^1$-local diffeomorphism on a compact Riemmanian manifold $M$
and assume that (robustly) every $f$-invariant ergodic probability measure $\mu$ satisfies the nonuniform specification property. Does the map $f$ satisfy the specification property?} \vspace{.2cm}

Another interesting topic is to relate specification properties with the presence of discontinuities 
in the system. Indeed, Buzzi~\cite{Buz97} proved that contrary to the
characterization due to Blokh~\cite{Bl83} for continuous interval maps
there exists a large class of topologically mixing but discontinuous maps of the interval
(including $\beta$-transformations) so that the set of parameters
for which the strong specification property holds although dense has
zero Lebesgue measure. In fact, it was communicated to us by Dan Thompson that there are 
beta-transformations that do not satisfy the non-uniform specification property for any 
full supported probability measure. So we pose the following question:

\vspace{.2cm}
{\bf Question 2:} \textit{Does the
set of parameters for which $\beta$-transformations satisfy the non-uniform specification property for any full supported invariant probability measure have positive Lebesgue measure?} 
 \vspace{.2cm}

Clearly this set contains  the set of parameters for which strong specification property holds, that has 
zero Lebesgue measure.  Note that an affirmative answer to the previous questions would be
a contribution for a better understanding of the non-uniform specification property would give a wider 
class of examples for which our results apply.

%%%%%%%%%%%%%%%%%%%%%%%
\subsubsection*{Large deviations}

Many recent contributions to the theory of large deviations in non-uniformly hyperbolic dynamical systems 
have been given. In fact, as discussed in the introduction, the existence of Markov towers allowed Melbourne, Nicol
~\cite{MN08} and Rey-Bellet, Young~\cite{RY08} to obtain large deviation principles with respect to the invariant
probability measure for the original dynamical system. 

More recently, Chung~\cite{Chu11} has also obtained large deviation principle with respect to the (not necessarily
invariant) Lebesgue measure for Markov tower maps induced by return time functions satisfying some technical conditions with similar flavor to our nonuniform specification property. Such results apply e.g. for a Markov
tower obtained from the Maneville-Pommeau with a rate function expressed by the pressure function computed using only expanding measures. The results in Example~\ref{ex.maneville} are a first step to prove that a large deviations principle hold for the original Maneville-Pommeau transformation. More generally, 
 
\vspace{.2cm}
{\bf Question 3:} \textit{Let $\mu$ be the unique equilibrium state for $f$ with respect to a potential $\phi$,
absolutely continuous with respect to a (not necessarily invariant) weak Gibbs measure, and whose first 
hyperbolic time map has exponential tail. Does there exists a large deviation principle with  respect to the 
weak Gibbs measure? If so, can the rate function be described using only thermodynamical quantities at
invariant expanding measures and the tail of the first hyperbolic time map?} \vspace{.2cm}

We expect that our results can extend to the partially hyperbolic setting. However, despite the fact that  
equilibrium states in a broad nonuniformly hyperbolic context are expected to be absolutely continuous 
with respect probability measures that exhibit some weak Gibbs property on Pesin local unstable leaves, 
the general thermodynamical formalism even for partially hyperbolic dynamical systems 
is far from being completely understood. 
Finally, let us mention that in the case of SRB measures some large deviations upper bounds were obtained in~\cite{AP06}.  Moreover,some  large deviations lower bounds have also been announced recently by Hirayama and
Sumi~\cite{HS10} for hyperbolic measures that satisfy a measure-theoretical transversality condition.

\section*{Acknowledgments} 
The author is grateful to V. Ara\'ujo, A. Castro and V. Pinheiro for their friendship and
advice and also to D. Thompson for interesting discussions on several weak specification properties.
The author is also grateful to the anonimous referee whose comments contributed to improve the manuscript.
This work was partially supported by CNPq-Brazil and FAPESB.

%%%%%%%%%%%%%%%%%%%%%%%%%%%%%%%%%%%%%%%%%%%%%%
\bibliographystyle{alpha}
\bibliography{bib}

\end{document}